\documentclass{amsart}

\usepackage{bbm}
\usepackage{amsfonts}
\usepackage{curves,multibox}
\usepackage{mathrsfs}
\usepackage{graphicx,verbatim}
\usepackage{amscd}
\usepackage{latexsym,tabularx}
\usepackage{amsmath}
\usepackage{amssymb}
\usepackage{amsthm}
\usepackage{color}
\usepackage{hyperref}
\newcommand{\GIM}{generalized intersection matrix}

\newcommand{\SGIM}{symmetrizable generalized intersection matrix}
\newcommand{\SIM}{symmetrizable intersection matrix}
\newcommand{\SIMs}{symmetrizable intersection matrices }

\newcommand{\SBF}{symmetric bilinear form}
\newtheorem{theorem}{Theorem}[section]

\newtheorem{lemma}[theorem]{Lemma}
\newtheorem{proposition}[theorem]{Proposition}
\newtheorem{corollary}[theorem]{Corollary}

\newcommand{\rank}{{\rm rank\,}}
\newcommand{\corank}{{\rm corank\,}}

\theoremstyle{definition}
\newtheorem{definition}[theorem]{Definition}

\theoremstyle{remark}
\newtheorem{remark}[theorem]{Remark}

\numberwithin{equation}{section}

\begin{document}
\title[Symmetrizable intersection matrices and their root systems] {Symmetrizable intersection matrices and their root systems}

\author{Liangang Peng}
\address{Department of Mathematics, Sichuan University, 610064 Chengdu, China}
\email{penglg@scu.edu.cn}
\thanks{Supported partially by the National Natural Science Foundation of China (No. 10325107),
and the 973 Program (No. 2006CB805905).}

\author{Mang Xu}
\address{Department of Mathematics, Southwest Jiaotong University, 610031 chengdu, China}
\email{xumang@home.swjtu.edu.cn}


\date{}


\begin{abstract}
   In this paper we study symmetrizable intersection matrices,
namely  generalized intersection matrices introduced by P. Slodowy
such that they are symmetrizable. Every such matrix can be naturally
associated with a root basis and a Weyl root system. Using $d$-fold
affinization matrices we give a classification, up to
braid-equivalence, for all  positive semi-definite symmetrizable
intersection matrices. We also give an explicit structure of the
Weyl root system for each $d$-fold affinization matrix in terms of
the root system of the corresponding Cartan matrix and some special
null roots.

\end{abstract}

\maketitle

\section{Introduction}
Generalized intersection matrices were introduced by P. Slodowy in
his study on isolated singularities (see \cite{Sl1} or \cite{Sl2}).
An integral matrix $A=(A_{ij})_{i, j\in I}$ with the finite index
set $I$ is called a generalized intersection matrix or a GIM if for
any $i, j\in I$, we have
$$\begin{array}{ccc}
A_{ii}=2, &             & \\
A_{ij}<0  & \mbox{iff}  & A_{ji}<0, \\
A_{ij}>0  & \mbox{iff}  & A_{ji}>0.
\end{array}$$
In this case, if $A$ is symmetrizable, then we call it a
symmetrizable intersection matrix or an SIM; and if $A$ is
symmetric, then we  call it an intersection matrix or an IM. Clearly
the notion of generalized intersection matrices is also a
generalization of the notion of generalized Cartan matrices.

For any generalized intersection matrix, one can associate it to a
root basis, as similar to what one dose for a generalized Cartan
matrix. Therefore a natural question is how to describe or classify
such kind of root bases.

When studying the elliptic singularities, K. Saito in \cite{Sai}
introduced the extended affine root systems, which are associated to
positive semi-definite bilinear forms such that the rank of the
radical of each such form is two. By using the diagrams with dotted
edges he obtained a classification for these extended affine root
systems.

B. Allison, N. Azam, S. Berman, Y. Gao and A. Pianzola in
\cite{AABGP} considered the extended affine root systems in general
form, namely these extended affine root systems associated to
arbitrary positive semi-definite bilinear forms without any
restriction on the rank of the radical of each such form. They used
the lattices or semi-lattices to obtain a classification for these
extended affine root systems.

In this paper we consider symmetrizable intersection matrices and
their root systems. The main method is applying
braid-transformations on root bases defined by Slodowy and diagrams
with dotted edges to our study. We state our main results as
follows.

Firstly, inspired by the Auslander-Platzeck-Reiten tilts or co-tilts
in representation theory of algebras we define the
APR-transformation on root bases as a composition of certain
braid-transformations. Using APR-transformations we  prove that any
positive semi-definite symmetrizable intersection matrix $A$ with
$\corank A\leq 1$ is braid-equivalent to a Cartan matrix or an
affine generalized Cartan matrix. In this case, as a consequence we
show that the intersection matrix Lie algebra of $A$ in sense of
Slodowy is isomorphic to a semi-simple or affine Kac-Moody Lie
algebra.

Secondly, given any Cartan Matrix, we define its (standard) $d$-fold
affinization matrix by adding some single roots or double roots of
this Cartan matrix. This generalizes the corresponding notion of S.
Berman and R. V. Moody \cite{BM} and that of  G. Benkart and E.
Zelmanov \cite{BZ}, where they only considered the case that single
roots were added. Furthermore, we prove that any  positive
semi-definite symmetrizable intersection matrix is braid-equivalent
to a $d$-fold affinization matrix. By introducing the type of a
$d$-fold affinization matrix and using standard $d$-fold
affinization matrices, we give a complete classification of all
positive semi-definite symmetrizable intersection matrices, up to
braid-equivalence.

Finally, we give an explicit structure of the Weyl root system for
each $d$-fold affinization matrix in terms of the root system of the
corresponding Cartan matrix and some special null roots.

We should mention that in this paper we do not deal with another
question how to define a suitable Lie algebra associated to a
generalized intersection matrix just like the Kac-Moody Lie algebra
associated to a generalized Cartan matrix. Usually this is a hard
question. There were many attempts to do it. For example, Slodowy in
\cite{Sl2} defined generalized intersection matrix Lie algebras and
intersection matrix Lie algebra; Berman and  Moody in \cite{BM} and
Benkart and Zelmanov in \cite{BZ} studied and classified the Lie
algebras graded by finite root systems; Allison, Azam, Berman, Gao
and  Pianzola in \cite{AABGP} considered the constructions of Lie
algebras to realize the extended affine root systems; Saito and
Yoshii in \cite{SY} defined simply-laced elliptic Lie algebras by
finite many generators and finite many relations depended on the
classification of elliptic root systems given by Saito in
\cite{Sai}; and so on. In \cite{XP}, we also made an attempt to
define the SIM Lie algebras for symmetrizable intersection matrices.

\section{SIM-Root Bases and APR-Transformations}

\subsection{GIM-Root Bases, SIM-Root Bases and Reflection
Transformations}

Let us first recall the definitions of generalized intersection
matrices, intersection matrices, GIM-root bases and
braid-equivalences given by P. Slodowy in \cite{Sl1}, \cite{Sl2}.
Then we consider symmetrizable generalized intersection matrices,
which are our emphasis.

\begin{definition}\label{GIM-IM}
{\rm Let $A=(A_{ij})_{i, j\in I}$ be an integral matrix, where $I$
is a finite set. If the integers $A_{ij}$ satisfies the following
properties: for all $i,j\in I$,
$$\begin{array}{ccc}
A_{ii}=2, &             & \\
A_{ij}<0  & \Leftrightarrow  & A_{ji}<0, \\
A_{ij}>0  & \Leftrightarrow  & A_{ji}>0,
\end{array}$$
then we call $A$ a {\it generalized intersection matrix}, or GIM for
short.

    If $A$ is symmetric, we call $A$ an {\it intersection
matrix}, or IM} for short.
\end{definition}

Note that any generalized Cartan matrix is a generalized
intersection matrix.

\begin{definition}\label{GIM-root}{  \rm
 Let $A=(A_{ij})_{i, j\in I}$ be a GIM.
 A {\it generalized intersection matrix root basis} of $A$, or a GIM-root basis for short,
 is a triple $(H, \nabla,  \Delta)$ consisting of\\
(1)\ a finite-dimensional $\mathbb{Q}$-vector space $H$;\\
(2)\ $\nabla =\{h_i | i\in I\}\subseteq H$;\\
(3)\ $\Delta=\{\alpha_i | i\in I \}\subseteq H^*={\rm
Hom}_{\mathbb{Q}}(H,{\mathbb{Q}})$, such that
$$\alpha_j (h_i)=A_{ij}\quad\mbox{for all }\ i, j\in I.$$

  We also call the generalized intersection matrix $A=(A_{ij})_{i, j\in I}$ is
 {\it the structural matrix} of the  GIM-root basis $(H, \nabla, \Delta)$. }
\end{definition}

  Let $A=(A_{ij})_{i, j\in I}$ be a GIM, and $( H, \nabla,
\Delta)$ be the corresponding GIM-root basis. The {\it reflection
transformations} are defined as follows: for every $\alpha\in
\Delta$, $$s_\alpha(h):=h-\alpha(h)h_\alpha,\quad h\in H.$$  The
contragredient action of $s_\alpha$ on $H^*$ is given by
$$s_\alpha(\gamma):=\gamma -\gamma (h_\alpha)\alpha,\quad\gamma\in
H^{\ast}.$$ Here if
 $\alpha=\alpha_i$, we denote $h_\alpha=h_i$.

  Two GIM-root bases $(H, \nabla, \Delta )$ and $(H, \nabla', \Delta')$
 are called {\it braid-equivalent} if they can be
transformed into each other by a sequence of transformations of the
form
$$(H, \nabla, \Delta )\mapsto\cdots\mapsto(H, \nabla_k, \Delta_k )\mapsto
(H, \nabla_{k+1}, \Delta_{k+1} )\mapsto\cdots\mapsto(H, \nabla',
\Delta' ),$$ where
$$\begin{array}{lcll}
\Delta_{k+1}&=&(\Delta_k\setminus\{\beta\})\cup\{s_{\alpha}(\beta)\},    & \\
            & &                                                          & \ \mbox{for\ some}\quad \alpha, \beta\in \Delta_k.\\
\nabla_{k+1}&=&(\nabla_k\setminus\{h_\beta\})\cup\{s_{\alpha}(h_\beta)\},&\\
\end{array}$$
It is easy to see that braid-equivalence is an equivalence relation.

  In the following we consider the symmetrizable case.

 Let $A=(A_{ij})_{i,j\in
I}$ be a {\it \SGIM}, that is, $A$ is a \GIM\, and there exists an
invertible diagonal matrix $D={\rm diag}(d_i)_{i\in I}$, where all
$d_i$'s are positive integers and $\gcd(d_i)_{i\in I} = 1$, such
that $DA$ is a symmetric matrix. In this case  $D$ is called the
{\it symmetrizer} of $A$. Denote $DA:= B=(B_{ij})_{i,j\in I}$. Let
$H$ be a $\mathbb{Q}$-vector space, and $\Delta :=\{\alpha_i | i\in
I \}\subseteq
 H$ be a linearly independent subset. Then there exists a \SBF
$$(-,\,-):\, H\times H \rightarrow \mathbb{Q}$$ such that
$$(\alpha_i,\,\alpha_j)=B_{ij}\quad\mbox{for}\ i, j\in I.$$
  It is easy to see that
$$A_{ij}=\frac{2(\alpha_i,\,\alpha_j)}{(\alpha_i,\,\alpha_i)}\quad\mbox{
for}\ i, j\in I$$ and $$d_i=(\alpha_i, \alpha_i)/2\quad\mbox{ for}\
i\in I.$$ We call $(H, \Delta)$ an {\it SIM-root basis} of the
\SGIM\, $A$, and $A$ is called the structural matrix of $(H,
\Delta)$. Note that since we have the above \SBF,  we can take
$\nabla=\{\frac{2\alpha_i}{(\alpha_i,\,\alpha_i)}| i\in I\}$.
Therefore we have no need to distinguish $H$ and its dual space
$H^*$. In this case, we have the following reflection
transformations: for every $\alpha_i\in\Delta$,
$$s_{\alpha_i}(h):=h- \frac{2(\alpha_i,\,
h)}{(\alpha_i,\,\alpha_i)}\alpha_i,\quad\ h\in H.$$ Particularly,
for all $\alpha_i,\ \alpha_j\in\Delta$,
$$s_{\alpha_i}(\alpha_j):=\alpha_j-A_{ij}\alpha_i.$$

    From now on we call a \SGIM\ just a {\it \SIM}, or SIM for short. We
always suppose $A=(A_{ij})_{i,j\in I}$ is an SIM. Accordingly, let
$H$, $\Delta$ and the \SBF\, $(-,\,-)$ be as above. Moreover, we
suppose $H$ is a $\mathbb{Q}$-vector space of dimension $|I|$. So
$\Delta$ is a basis of $H$. In this case we say $\Delta$ is an
SIM-root basis of the \SIM\,$A$, and $A$ is the structural matrix of
$\Delta$.

If $A=(A_{ij})_{i, j\in I}$ is a (symmetrizable) generalized Cartan
matrix, we call the corresponding SIM-root basis or GIM-root basis a
{\it GCM-root basis}.

 We fix the following notations and notions.

 Let ${\mathcal W}\subseteq GL(H)$ be the subgroup of $GL(H)$ generated by $s_{\alpha_i}$, $\alpha_i\in
   \Delta$. We call $\mathcal{W}$ the {\it Weyl group} of the SIM-root basis
   $\Delta$, and
$$R^{W}:=\{s(\alpha)\,|\, s\in {\mathcal W}, \alpha\in \Delta
\}$$ the {\it Weyl root system} of $\Delta$. Set
$$\Gamma:=\sum\limits_{i\in I}\mathbf{Z}\alpha_i.$$
The lattice $\Gamma$ is called the {\it root lattice} of $\Delta$.

If an element $\alpha\in H$ is non-isotropic, that is
$(\alpha,\,\alpha)\not=0$, then we define the dual $\alpha^\vee$ of
$\alpha$ and the reflection $s_\alpha\in GL(H)$ as follows:
$$\alpha^\vee:=\frac{2}{(\alpha,\,\alpha)}\alpha,$$
$$s_\alpha(\beta):=\beta-(\beta,\,\alpha^\vee)\alpha,
\quad\mbox{for all }\,\beta\in H.$$

It is easy to check that
$$\alpha^{\vee\vee}=\alpha,$$
$$s_\alpha=s_{\alpha^{\vee}},\quad s_\alpha^2=id_H,$$
$$(s_\alpha(\beta),\,s_\alpha(\gamma))=(\beta,\,\gamma),\quad s_\alpha(\beta)^\vee=s_\alpha(\beta^\vee).$$

\subsection{Braid-transformations and APR-transformations}

Let $\Delta$ be an SIM-root basis with a \SIM\ $A=(A_{ij})_{i, j\in
I}$ as its structural matrix.

For all $a, b\in I$, we define an `operator' $\tau_{a, b}$ on the
root basis $\Delta$  as follows:
$$\tau_{a, b}\Delta=\Delta'$$ where
$\Delta'=\{\,\alpha'_i\,|\,i\in I \}$ is defined as:
$$\begin{array}{ll}
\alpha'_i=\alpha_i &\mbox{ for all } b\neq i\in I,  \\
\alpha'_b=s_{\alpha_a}(\alpha_b)=\alpha_b-A_{ab}\alpha_a. &
\end{array}$$
The structural matrix of  $\Delta'$ is  $A'=(A'_{ij})_{i, j\in I}$,
where
$$\begin{array}{ll}
A'_{ij}=A_{ij} &\quad\mbox{for} \ \ i\neq b,j\neq b, \\
A'_{ib}=A_{ib}-A_{ia}A_{ab} &\quad\mbox{for} \ \ i\neq b,\\
A'_{bj}=A_{bj}-A_{ba}A_{aj} &\quad\mbox{for} \ \ j\neq b,\\
A'_{bb}=2. &
\end{array}$$
 It is easy to see that $A'$ is also a \SIM\ and its symmetrizer is the same as that of $A$.
 Therefore, $\Delta'$ is also an SIM-root basis.
 We call $\tau_{a, b}$ a {\it braid-transformation}, and write $\tau_{a, b}A:=A'$.

Two SIM-root bases $\Delta$ and $\Delta'$ are called {\it
braid-equivalent}, if there exist $a_1, b_1, \cdots,$ $ a_t, b_t \in
I$ such that
$$\tau_{a_t, b_t}\cdots \tau_{a_1, b_1}\Delta= \Delta'.$$ Obviously, the braid-transformation $\tau_{a,
b}$ is invertible, and its inverse is also $\tau_{a, b}$. Therefore
the braid-equivalence relation is  an equivalence relation. In fact,
the definition of the braid-equivalences of SIM-root bases given
here is coincident with that of the corresponding braid-equivalences
of GIM-root bases.

If the structural matrices of the above two SIM-root bases are
respectively $A=(A_{ij})_{i, j\in I}$ and $A'=(A'_{ij})_{i, j\in
I}$, we  also say $A$ and $A'$  braid-equivalent. It is easy to show
that a braid-transformation preserves the property of being
symmetric, symmetrizable or indecomposable as well as the rank of a
matrix.

\begin{remark}
{\rm In our definition, the structural matrix $A=(A_{ij})_{i, j\in
I}$ of an SIM-root basis or a GIM-root basis is independent of the
order of the finite index set $I$. Therefore two \SIMs
 $A=(A_{ij})_{i, j\in I}$ and $A'=(A'_{ij})_{i, j\in I}$ of the same
order of the finite index set are braid-equivalent if and only if
$A$ can be transformed into $A'$ by a sequence of `operators' of
this form $\tau_{a, b}$($a, b\in I$) and permutations of elements in
the index set $I$.}
\end{remark}

    More generally, for $a\in I$ and a nonempty subset $N$ of $I$, we define the operator $\tau_{a,\,N}$
as follows:
   $$\tau_{a,\,N}\Delta=\Delta',$$ where $\Delta'=\{\alpha'_i\;|\;i\in
I\}$ with
$$\begin{array}{ll}
 \alpha'_i=\alpha_i, & \quad\mbox{for all } \ i\in
I\setminus N;\\
\alpha'_i=s_{\alpha_a}(\alpha_i), & \quad \mbox{for all } \ i\in N.\\
\end{array}$$
The operator $\tau_{a,\,N}$ has the following property:

\begin{proposition}\label{braid}
Let $N=\{ b_1, b_2, \cdots, b_t \} \subseteq I$. Then
$$\tau_{a,\,N}\Delta=\tau_{a,\,b_t}\cdots\tau_{a,\,b_1}\Delta.$$
\end{proposition}
\begin{proof}
We prove it by induction. For $t=1$, it is clear. For $t>1$, put
$N'=N\setminus\{b_t\}$ and suppose that
$\tau_{a,\,N'}\Delta=\tau_{a,\,b_{t-1}}\cdots\tau_{a,\,b_1}\Delta$.
Set $\tau_{a,\,N'}\Delta=\Delta'$. Then $\alpha'_i=\alpha_i$ if
$i\in I\setminus N'$, and $\alpha'_i=s_{\alpha_a}(\alpha_i)$ if
$i\in N'$. Set $\tau_{a,\,b_t}\Delta'=\Delta''$. Then
$\alpha''_i=\alpha'_i$ if $i\neq b_t$, and
$\alpha''_{b_t}=s_{\alpha_a'}(\alpha'_{b_t})
=\alpha'_{b_t}-(\alpha'^{\vee}_a,\,\alpha'_{b_t})\alpha'_a$. Thus if
$i\in I\setminus N$, then $\alpha''_i=\alpha_i$. If $i\in N'$, then
$\alpha''_i=s_{\alpha_a}(\alpha_i)$. So, when $a\in I\setminus N'$,
we have
$\alpha''_{b_t}=\alpha_{b_t}-(\alpha_a^{\vee},\,\alpha_{b_t})\alpha_a=s_{\alpha_a}(\alpha_{b_t})$.
And when $a\in N'$, we have $\alpha''_{b_t}
=\alpha_{b_t}-(-\alpha^{\vee}_a,\,\alpha_{b_t})(-\alpha_a)=s_{\alpha_a}(\alpha_{b_t})$.
Therefore we have
$\tau_{a,\,N}\Delta=\Delta''=\tau_{a,\,b_t}\Delta'=\tau_{
a,\,b_t}\cdots\tau_{a,\,b_1}\Delta$.
\end{proof}

\begin{definition}\label{APR}
{  \rm  Let $A=(A_{ij})_{i,j\in I}$ be an SIM, and $\Delta$ be the
SIM-root basis of $A$. Given $a\in I$, let $N(a)=\{ a\}\cup \{ b\in
I |\ A_{ab}<0\}$. We call $\tau_a:=\tau_{a,\,N(a)}$ an {\it
Auslander-Platzeck-Reiten transformation} of the root basis $\Delta$
at the point $a$, or an APR-transformation at $a$ for brevity.}
\end{definition}

\begin{remark}
{\rm  The notion of Auslander-Platzeck-Reiten transformations is
inspired by the notion of APR-tilts and APR-co-tilts given by M.
Auslander, M. I. Platzeck and I. Reiten in \cite{APR}, which  plays
 an important role in the representation theory of algebras.  From
 the view of tilting theory in the representation theory of
 algebras, L. Peng in \cite{Pe} gave an explicit explanation of
the relation between APR-transformations given here and APR-tilts as
well as APR-co-tilts.}
\end{remark}

Proposition~\ref{braid} says that an APR-transformation is a
composite of certain braid-transformations. In the following we
shall produce a formula to compute the structural matrix of an
APR-transformation of an SIM-root basis. This formula is very useful
for our computations.

\begin{proposition}\label{APReq}
Let $\tau_a$ be an APR-transformation. Denote $\Delta^a
:=\tau_{a}\Delta $ with structural matrix $A^a\ (:=\tau_{a}A)$,
where $A^a=(A_{ij}^a)_{i,j\in I}$ is the matrix such that for
 all $i,j\in I$, $A^a_{ij}=(\alpha^a_j ,\ (\alpha_i^a)^\vee )$. Then we have
$$\begin{array}{ll}
A^a_{ij}=A_{ij} &\quad\mbox{for} \ \ i,j\not\in N(a)\ \ \mbox{or}\ \
i,j\in N(a), \\
A^a_{ij}=A_{ij}-A_{ia}A_{aj} &\quad \mbox{otherwise}.
\end{array}$$
\end{proposition}

\begin{proof}
 By definition, for $i\not\in N(a)$,  we have $\alpha^a_i=\alpha_i$; for $i\in
N(a)$, we have $\alpha^a_i=s_{\alpha_a}(\alpha_i)=\alpha_i-
A_{ai}\alpha_a$. So when $i,j\not\in N(a)$ we have
$$A^a_{ij}=(\alpha^a_j ,\ (\alpha_i^a)^\vee )=(\alpha_j ,\
\alpha_i^\vee )=A_{ij}.$$ When $i,j\in N(a)$ we have
$$A^a_{ij}=(s_{\alpha_a}(\alpha_j) ,\ s_{\alpha_a}(\alpha_i)^\vee
)=(\alpha_j ,\ \alpha_i^\vee )=A_{ij}.$$ When $i\in N(a)$ and
$j\not\in N(a)$ we have
\begin{eqnarray*}
A^a_{ij} & = & (\alpha_j ,\
s_{\alpha_a}(\alpha_i)^\vee)=\frac{2(\alpha_j, \
\alpha_i-A_{ai}\alpha_a)}{(s_{\alpha_a}(\alpha_i), \
s_{\alpha_a}(\alpha_i) )}\\
& = & A_{ij}-\frac{2A_{ai}(\alpha_j ,\ \alpha_a )}{(\alpha_i , \
\alpha_i )}=A_{ij}-\frac{4(\alpha_i ,\ \alpha_a )(\alpha_j ,\
\alpha_a )}{(\alpha_a ,\ \alpha_a )(\alpha_i , \ \alpha_i
)}=A_{ij}-A_{ia}A_{aj}.
\end{eqnarray*}
When $i\not\in N(a)$ and $ j\in N(a)$, we have
$$A^a_{ij}=(s_{\alpha_a}(\alpha_j) ,\ \alpha_i^\vee
)=\frac{2(\alpha_j-A_{aj}\alpha_a, \ \alpha_i )}{(\alpha_i, \
\alpha_i )}
 = A_{ij}-\frac{2A_{aj}(\alpha_a ,\ \alpha_i)}{(\alpha_i , \
\alpha_i )}= A_{ij}-A_{ia}A_{aj}.$$
\end{proof}

By the definition of $N(a)$ and the above formula, we have the
following immediate consequence.

\begin{corollary}
If $A$ is a symmetrizable generalized Cartan matrix, then for any
APR-transformation $\tau_a$, we have $\tau_a A =A$. \hfill$\square$
\end{corollary}

\subsection{SIM's and Dynkin Diagrams}

For any \SIM\ $A=(A_{ij})_{i, j\in I}$ we define a diagram $D(A)$
associated to $A$, called the Dynkin diagram of $A$, as follows.

Let $I$ be the vertex set. If $a_{ij}\neq 0$ and $a_{ij}\neq
a_{ji}$, the vertices $i$ and $j$ are connected by a solid edge if
$a_{ij}<0$ or by a dotted edge if $a_{ij}>0$, and these edges are
equipped with an ordered pair of integers $(|a_{ij}|, |a_{ji}|)$. If
$a_{ij}=a_{ji}$, the vertices $i$ and $j$ are connected by
$(-a_{ij})$ solid edges if $a_{ij}\leq 0$ or by $a_{ij}$ dotted
edges if $a_{ij}>0$. Here, our definition is slightly different from
that given by Slodowy in \cite{Sl1} and \cite{Sl2}. Our edges are
valued, while those defined by Slodowy are not valued.

    For two SIM's $A$ and $A'$, $D(A)$ and $D(A')$ are said to be braid-equivalent
if their SIM-root bases are braid-equivalent. In addition, if
$\tau_{a,\,b}A=A'$, then we also write $\tau_{a,\,b}D(A)=D(A')$.

    A diagram without dotted edges is called a {\it solid diagram}. It is
clear that $D(A)$ is a solid diagram if and only if $A$ is a
generalized Cartan matrix. In this case $D(A)$ is invariant under
any APR-transformation.

    For the sake of convenience, we label
all Dynkin diagrams of positive-definite GCM's (that is, classic
Dynkin diagrams) as follows:

$A$ type:
$$\unitlength=1cm
\begin{picture}(4,0.6)
\put(0,0.27){$a_t$} \put(1.4,0.26){$\cdots$} \put(2.8,0.27){$a_2$}
\put(4.3,0.27){$a_1$} \put(0.4,0.35){\line(1,0){0.7}}
\put(2,0.35){\line(1,0){0.7}} \put(3.3,0.35){\line(1,0){0.7}}
\end{picture}$$

$B$ type:
$$\unitlength=1cm
\begin{picture}(6.5,0.6)
\put(0,0.27){$a_t$} \put(0.5,0.35){\line(1,0){0.7}}
\put(1.4,0.27){$a_{t-1}$} \put(2.2,0.35){\line(1,0){0.7}}
\put(3.1,0.27){$\cdots$} \put(3.8,0.35){\line(1,0){0.7}}
\put(4.7,0.27){$a_0$} \put(5.2,0.35){\line(1,0){0.7}}
\put(5.55,0.5){\makebox(0,0)[c]{\scriptsize $(1,2)$}}
\put(6.1,0.27){$b$}
\end{picture}$$

$C$ type:
$$\unitlength=1cm
\begin{picture}(6.5,0.6)
\put(0,0.27){$a_t$} \put(0.5,0.35){\line(1,0){0.7}}
\put(1.4,0.27){$a_{t-1}$} \put(2.2,0.35){\line(1,0){0.7}}
\put(3.1,0.27){$\cdots$} \put(3.8,0.35){\line(1,0){0.7}}
\put(4.7,0.27){$a_0$} \put(5.2,0.35){\line(1,0){0.7}}
\put(5.55,0.5){\makebox(0,0)[c]{\scriptsize $(2,1)$}}
\put(6.1,0.27){$b$}
\end{picture}$$

$D$ type:
$$\unitlength=1cm
\begin{picture}(6,1.5)
\put(0,0.27){$a_t$} \put(0.5,0.35){\line(1,0){0.7}}
\put(1.4,0.27){$\cdots$} \put(2.1,0.35){\line(1,0){0.7}}
\put(3,0.27){$a_1$} \put(3.5,0.35){\line(1,0){0.7}}
\put(4.4,0.27){$a_0$} \put(4.9,0.35){\line(1,0){0.7}}
\put(5.8,0.27){$b$} \put(4.45,1.5){$c$}
\put(4.5,1.3){\line(0,-1){0.7}}
\end{picture}$$

$E_6$ type:
$$\unitlength=1cm
\begin{picture}(6,1.5)
\put(0,0.27){$a_2$} \put(0.5,0.35){\line(1,0){0.7}}
\put(1.4,0.27){$a_1$} \put(1.9,0.35){\line(1,0){0.7}}
\put(2.8,0.27){$a_0$} \put(3.3,0.35){\line(1,0){0.7}}
\put(4.2,0.27){$b_1$} \put(4.7,0.35){\line(1,0){0.7}}
\put(5.6,0.2){$b_2$} \put(2.9,1.5){$c$}
\put(3,1.3){\line(0,-1){0.7}}
\end{picture}$$

$E_7$ type:
$$\unitlength=1cm
\begin{picture}(6,1.5)
\put(0,0.27){$a_3$} \put(0.5,0.35){\line(1,0){0.7}}
\put(1.4,0.27){$a_2$} \put(1.9,0.35){\line(1,0){0.7}}
\put(2.8,0.27){$a_1$} \put(3.3,0.35){\line(1,0){0.7}}
\put(4.2,0.27){$a_0$} \put(4.7,0.35){\line(1,0){0.7}}
\put(5.6,0.27){$b_1$} \put(6.1,0.35){\line(1,0){0.7}}
\put(7.0,0.27){$b_2$} \put(4.3,1.5){$c$}
\put(4.4,1.3){\line(0,-1){0.7}}
\end{picture}$$

$E_8$ type:
$$\unitlength=1cm
\begin{picture}(7,1.5)
\put(0,0.27){$a_4$} \put(0.5,0.35){\line(1,0){0.7}}
\put(1.4,0.27){$a_3$} \put(1.9,0.35){\line(1,0){0.7}}
\put(2.8,0.27){$a_2$} \put(3.3,0.35){\line(1,0){0.7}}
\put(4.2,0.27){$a_1$} \put(4.7,0.35){\line(1,0){0.7}}
\put(5.6,0.27){$a_0$} \put(6.1,0.35){\line(1,0){0.7}}
\put(7,0.27){$b_1$} \put(7.5,0.35){\line(1,0){0.7}}
\put(8.4,0.27){$b_2$} \put(5.7,1.5){$c$}
\put(5.8,1.3){\line(0,-1){0.7}}
\end{picture}$$

$F_4$ type:
$$\unitlength=1cm
\begin{picture}(5,0.6)
\put(0,0.27){$a_1$} \put(0.5,0.35){\line(1,0){0.7}}
\put(1.4,0.27){$a_0$} \put(2,0.35){\line(1,0){0.7}}
\put(2.35,0.5){\makebox(0,0)[c]{\scriptsize $(1,2)$}}
\put(2.9,0.27){$b_1$} \put(3.45,0.35){\line(1,0){0.7}}
\put(4.3,0.27){$b_2$}
\end{picture}$$

$G_2$ type:
$$\unitlength=1cm
\begin{picture}(2,0.6)
\put(0,0.27){$a_0$} \put(0.45,0.35){\line(1,0){0.7}}
\put(0.8,0.5){\makebox(0,0)[c]{\scriptsize $(1,3)$}}
\put(1.4,0.27){$b$}
\end{picture}$$

    Let $a_t$ be the left end point in each Dynkin diagram as above, for example,
    $a_t=a_2$ in $E_6$ case, similarly for others. Then we have

\begin{lemma} \label{SIM-GCM}
Let $A=(A_{ij})_{i,j\in I}$ be an SIM such that $A$ is positive
semi-definite and $\corank A\leq 1$. If there is a $d\in I$ such
that the full subdiagram $D(A)\setminus \{d\}$ of $D(A)$ is a Dynkin
diagram of type $A_l, B_l, C_l, D_l, E_6, E_7, E_8, F_4, G_2$ with
$A_{di}\leq 0$ for all $i\in I\setminus \{d,a_t \}$, then the
SIM-root basis $\Delta$ is braid-equivalent to a GCM-root basis by
successive APR-transformations at some vertices in
$I\setminus\{d\}$.
\end{lemma}

\begin{proof}
 If $A_{da_t}\leq 0$, then $A$ itself is a GCM. So
we can assume that $A_{da_t}>0$. Note that $A$ is positive
semi-definite. Thus we have $A_{da_t}A_{a_t d}\leq 4$ and $h\leq 4$,
where $h=|\sum\limits_{i\in I\setminus\{d,a_t\}}A_{di}|$.

    In the following diagrams, the valuation on the vertical edge connecting
the top vertex $i$ and the bottom vertex $j$ is denoted $(a_{ij},
a_{ji})$. In other cases, $(a_{ij}, a_{ji})$ means the valuation on
the edge linking the left vertex $i$ and the right  vertex $j$.

    Since $A_{da_t}A_{a_t d}\leq 4$, we have $A_{da_t}\leq 4$. So we verify the statement case by case,
according to the value of $A_{da_t}$ being equal to $4,3,2$ or $1$,
respectively. We will not list some obviously impossible cases,
taking into account the condition that our matrix is symmetrizable
and positive semi-definite.\\

    {\it Case 1}. Assume that $A_{da_t}=4$. Then we have $A_{a_td}=1$. There are
    three sub-cases.

    (1.1) Let $D(A)\setminus\{d\}$ be of type $A_l$.

    It is clear that $D(A)\setminus\{d\}$
must be of type $A_1$. Thus $\tau_{a_1}\Delta$ is a GCM-root basis
of type $A^{(2)}_2$.

    (1.2) Assume that $D(A)\setminus\{d\}$ is of type $B_l$.

    It is easy to see there is only one case:\\
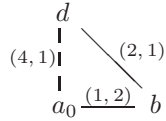
$$\unitlength=1cm
$$
The diagram of $\tau_b\tau_{a_0}\Delta$ is
$$\unitlength=1cm
%
$$

    (1.3) $D(A)\setminus\{d\}$ can not be of type $C_l(l>2) , D_l, E_6, E_7, E_8, F_4, G_2$.\\

    {\it Case 2}.  Assume that $A_{da_t}=3$. Then we have $A_{a_td}=1$. There
    are also three sub-cases.

    (2.1) Let $D(A)\setminus\{d\}$ be of type $A_l$.

    It is easy to see that $D(A)\setminus\{d\}$ must be of type $A_1$ or $A_2$.

    If it is $A_1$, $\tau_{a_1}\Delta$ is a GCM-root basis of type
$G_2$.

    If it is $A_2$, there are only two possibilities. One is
$$\unitlength=1cm
%
$$
The diagram of $\tau_{a_1}\tau_{a_2}\Delta$ is
$$\unitlength=1cm
%
$$
 Another one is
$$\unitlength=1cm
%
$$
The diagram of $\tau_{a_2}\Delta$ is
$$\unitlength=1cm
%
$$

    (2.2)  $D(A)\setminus\{d\}$ can not be of type $B_l, C_l, D_l, E_6, E_7, E_8, F_4$.

    (2.3)  Suppose $D(A)\setminus\{d\}$ is of type $G_2$. There are only two
    possibilities. The first one is
$$\unitlength=1cm
%
$$
The diagram of $\tau_b\tau_{a_0}\tau_b\tau_{a_0}\Delta$ is
$$\unitlength=1cm
%
$$
The second one is
$$\unitlength=1cm
%
$$
The diagram of $\tau_b\tau_{a_0}\Delta$ is
$$\unitlength=1cm
%
$$

    {\it Case 3}. Assume that $A_{da_t}=2$ and $A_{a_td}=1$. Then we
    have the following six sub-cases.

    (3.1) Suppose $D(A)\setminus\{d\}$ is of type $A_l$. There are two
    possibilities, according to the value of $h$.

    (3.1.1) If $h=0$, the diagram of $\Delta$ is
$$\unitlength=1cm
%
$$
The diagram of $\tau_{a_1}\tau_{a_2}\cdots\tau_{a_t}\Delta$ is
$$\unitlength=1cm
%
$$

    (3.1.2) If $h=2$, the diagram of $\Delta$ is
$$\unitlength=1cm
%
$$
The diagram of $\tau_{a_{i+1}}\tau_{a_{i+2}}\cdots\tau_{a_t}\Delta$
is
$$\unitlength=1cm
%
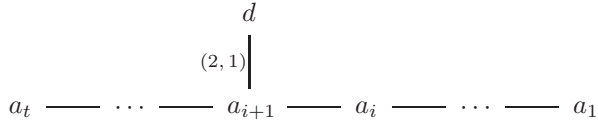
$$

    (3.2)  Assume that $D(A)\setminus\{d\}$ is of type $B_l$. Then there
    are three possibilities.

    (3.2.1) If $h=0$, the diagram of $\Delta$ is
$$\unitlength=1cm
%
$$
The diagram of
$\tau_{a_{t}}\cdots\tau_{a_0}\tau_b\tau_{a_{0}}\cdots\tau_{a_t}\Delta$
is
$$\unitlength=1cm
%
$$

    (3.2.2) If $h=1$, the diagram of $\Delta$ is only one
    possibility:
$$\unitlength=1cm
%
$$
The diagram of $\tau_b\tau_{a_0}\cdots\tau_{a_t}\Delta$ is
$$\unitlength=1cm
%
$$

    (3.2.3) If $h=2$, then there are two possibilities. One is
$$\unitlength=1cm
%
$$
The diagram of $\tau_{a_{t}}\Delta$ is
$$\unitlength=1cm
%
$$
Another one is
$$\unitlength=1cm
%
$$
The diagram of $\tau_{a_{0}}\Delta$ is
$$\unitlength=1cm
%
$$

    (3.3) Suppose $D(A)\setminus\{d\}$ is of type $C_l$. Similarly, according
    to the value of $h$, we have the following cases.

    (3.3.1) If $h=0$, the diagram of $\Delta$ is
$$\unitlength=1cm
%
$$
The diagram of
$\tau_{a_{t}}\cdots\tau_{a_0}\tau_b\tau_{a_{0}}\cdots\tau_{a_t}\Delta$
is
$$\unitlength=1cm
%
$$

    (3.3.2) If $h=2$, the diagram of $\Delta$ is only one
    possibility:
$$\unitlength=1cm
%
$$
The diagram of $\tau_{a_{t}}\Delta$ is
$$\unitlength=1cm
%
$$

    (3.3.3) If $h=4$, the diagram of $\Delta$ is only one
    possibility:
$$\unitlength=1cm
%
$$
The diagram of $\tau_{a_{0}}\Delta$ is
$$\unitlength=1cm
%
$$

    (3.4) Assume that $D(A)\setminus\{d\}$ is of type $D_l$. We have the
    following cases.

    (3.4.1) If $h=0$, the diagram of $\Delta$ is
$$\unitlength=1cm
%
$$
The diagram of
$\tau_{a_t}\cdots\tau_{a_0}\tau_c\tau_b\tau_{a_0}\cdots\tau_{a_t}\Delta$
is
$$\unitlength=1cm
%
$$

    (3.4.2) If $h=2$, there are two cases. One is
$$\unitlength=1cm
%
$$
The diagram of $\tau_{a_t}\Delta$ is
$$\unitlength=1cm
%
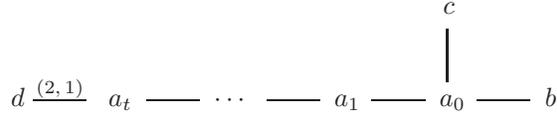
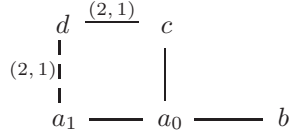
$$
Another one is
$$\unitlength=1cm
%
$$
The diagram of $\tau_{b}\tau_{a_0}\tau_{a_1}\Delta$ is
$$\unitlength=1cm
%
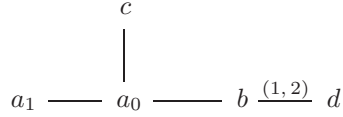
$$

    (3.5)  $D(A)\setminus\{d\}$ can not be of type $E_l (l=6,7,8)$ or $G_2$.

    (3.6) Suppose $D(A)\setminus\{d\}$ is of type $F_4$. Then $h$  must be equal to 1
and there are only two possibilities. One is
$$\unitlength=1cm
%
$$
The diagram of $\tau_{b_2}\tau_{b_1}\tau_{a_0}\tau_{a_1}\Delta$ is
$$\unitlength=1cm
%
$$
The diagram of $\tau_{b_2}\tau_{b_1}\tau_{a_0}\tau_{a_1}\Delta$ is
$$\unitlength=1cm
%
$$

    {\it Case 4}. Assume that $A_{da_t}=2$ and $A_{a_td}=2$. Then we
    have the following sub-cases.

    (4.1) Assume that $D(A)\setminus\{d\}$ is of type $A_l$.

    If it is $A_1$, $\tau_{a_1}\Delta$ is a GCM-root
basis of type $A^{(1)}_1$.

   If $l>1$, then $h$ must be equal to 1. There is only one possibility

$$\unitlength=1cm
%
$$
The diagram of $\tau_{a_1}\cdots\tau_{a_t}\Delta$ is
$$\unitlength=1cm
%
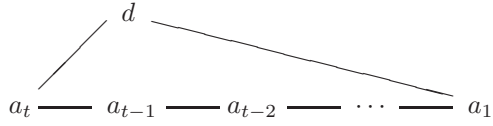
$$

    (4.2) Assume that $D(A)\setminus\{d\}$ is of type $B_l$. Then $h$ must be equal to 1
and there are two possibilities. One is
$$\unitlength=1cm
%
$$
The diagram of
$\tau_{a_t}\cdots\tau_{a_0}\tau_{b}\tau_{a_0}\cdots\tau_{a_t}\Delta$
is
$$\unitlength=1cm
%
$$
The diagram of $\tau_b\tau_{a_{0}}\Delta$ is
$$\unitlength=1cm
%
$$

    (4.3) Suppose $D(A)\setminus\{d\}$ is of type $C_l$.

    (4.3.1) If $h=1$, the diagram of $\Delta$ is only one
    possibility:
$$\unitlength=1cm
%
$$
The diagram of
$\tau_{a_t}\cdots\tau_{a_0}\tau_b\tau_{a_0}\cdots\tau_{a_t}\Delta$
is
$$\unitlength=1cm
%
$$
The diagram of $\tau_{b}\tau_{a_{0}}\Delta$ is

$$\unitlength=1cm
%
$$

    (4.4) Assume that $D(A)\setminus\{d\}$ is of type $D_l$. Then the diagram of $\Delta$ is only one
    possibility:
$$\unitlength=1cm
%
$$
The diagram of
$\tau_{a_{t-1}}\cdots\tau_{a_0}\tau_c\tau_b\tau_{a_0}\cdots\tau_{a_t}\Delta$
is
$$\unitlength=1cm
%
$$

    (4.5) Suppose $D(A)\setminus\{d\}$ is of type $E_l$. The diagram of $\Delta$ is only one
    possibility:
$$\unitlength=1cm
%
$$
For $t=2$, the diagram of
$\tau_c\tau_{a_0}\tau_{b_1}\tau_{a_1}\tau_{a_0}\tau_c\tau_{b_2}\tau_{b_1}\tau_{a_0}
\tau_{a_1}\cdots\tau_{a_t}\Delta$ is
$$\unitlength=1cm
%
$$
For $t=3$, the diagram of $\tau_{b_2}\tau_{b_1}\tau_{a_0}\tau_{a_1}
\tau_c\tau_{a_0}\tau_{a_2}\tau_{a_1}\tau_{b_1}\tau_{a_0}
\tau_c\tau_{b_2}\tau_{b_1}\tau_{a_0}\tau_{a_1}\cdots\tau_{a_t}\Delta$
is
$$\unitlength=1cm
%
$$
For $t=4$, the diagram of
$\tau_{a_4}\tau_{a_3}\tau_{a_2}\tau_{a_1}\tau_{a_0}\tau_{b_1}\tau_c\tau_{a_0}\tau_{a_1}
\tau_{b_2}\tau_{b_1}\tau_{a_0}\tau_{a_2}\tau_{a_1}\tau_c\tau_{a_0}\tau_{a_3}\tau_{a_2}
\tau_{a_1}\\
\tau_{b_1}\tau_{a_0}\tau_c\tau_{b_2}\tau_{b_1}\tau_{a_0}
\tau_{a_1}\cdots\tau_{a_t}\Delta$ is
$$\unitlength=1cm
%
$$

    (4.6) Assume that $D(A)\setminus\{d\}$ is of type $F_4$. The diagram of $\Delta$ is only one
    possibility:
$$\unitlength=1cm
%
$$
The diagram of
$\tau_{a_1}\tau_{a_0}\tau_{b_1}\tau_{b_2}\tau_{a_0}\tau_{b_1}\tau_{a_0}\tau_{a_1}\Delta$
is
$$\unitlength=1cm
%
$$

    (4.7) Assume that $D(A)\setminus\{d\}$ is of type $G_2$. And the diagram of $\Delta$ is only one
    possibility:
$$\unitlength=1cm
%
$$
The diagram of $\tau_{a_0}\tau_b\tau_{a_{0}}\Delta$ is
$$\unitlength=1cm
%
$$

    {\it Case 5}. Assume that $A_{d a_t}=1$ and $A_{a_t d}=1$. Then we
    have the following sub-cases.

    (5.1) Assume that $D(A)\setminus\{d\}$ is of type $A_l$. Then we have
    the following three cases.

    (5.1.1) If $h=0$, the diagram of $\Delta$ is
$$\unitlength=1cm
%
$$
The diagram of $\tau_{a_1}\cdots\tau_{a_{t-1}}\tau_{a_t}\Delta$ is
$$\unitlength=1cm
%
$$

    (5.1.2) If $h=1$, the diagram of $\Delta$ is
$$\unitlength=1cm
%
$$
The diagram of $\tau_{a_{i+1}}\tau_{a_{i+2}}\cdots\tau_{a_t}\Delta$
is
$$\unitlength=1cm
%
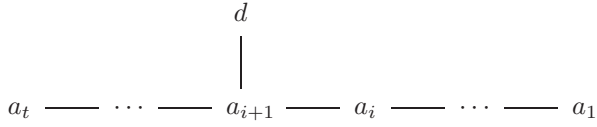
$$

    (5.1.3) If $h=2$, there are two cases. One is
$$\unitlength=1cm
%
$$
The diagram of $\tau_{a_t}\Delta$ is
$$\unitlength=1cm
%
$$
The diagram of $\tau_{a_2}\Delta$ is
$$\unitlength=1cm
%
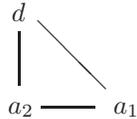
$$

    (5.2) Assume that $D(A)\setminus\{d\}$ is of type $B_l$. Then
    similarly we have the following cases.

    (5.2.1) If $h=0$, the diagram of $\Delta$ is
$$\unitlength=1cm
%
$$
The diagram of
$\tau_{a_{t}}\cdots\tau_{a_0}\tau_b\tau_{a_{0}}\cdots\tau_{a_t}\Delta$
is
$$\unitlength=1cm
%
$$

    (5.2.2) If $h=1$, then there are two possibilities. One is
$$\unitlength=1cm
%
$$
The diagram of $\tau_{a_{i+1}}\tau_{a_{i+2}}\cdots\tau_{a_t}\Delta$
is
$$\unitlength=1cm
%
$$
The diagram of $\tau_{a_0}\cdots\tau_{a_t}\Delta$ is
$$\unitlength=1cm
%
$$

    (5.3) Assume that $D(A)\setminus\{d\}$ is of type $C_l$. Then we have the following cases.

    (5.3.1) If $h=0$, the diagram of $\Delta$ is
$$\unitlength=1cm
%
$$
The diagram of $\tau_{a_0}\cdots\tau_{a_t}\Delta$ is
$$\unitlength=1cm
%
$$

    (5.4) Suppose $D(A)\setminus\{d\}$ is of type $D_l$. As before we have the following cases
    according to the value of $h$.

    (5.4.1) If $h=0$, the diagram of $\Delta$ is
$$\unitlength=1cm
%
$$
The diagram of
$\tau_{a_t}\cdots\tau_{a_1}\tau_{a_0}\tau_b\tau_c\tau_{a_0}\tau_{a_1}\cdots\tau_{a_t}\Delta$
is
$$\unitlength=1cm
%
$$

    (5.4.2) If $h=1$, then there are two possibilities. One is
$$\unitlength=1cm
%
$$
The diagram of
$\tau_b\tau_{a_0}\tau_{a_1}\tau_{a_2}\cdots\tau_{a_t}\Delta$ is
$$\unitlength=1cm
%
$$

    (5.5) Suppose $D(A)\setminus\{d\}$ is of type $E_l$. We have the following cases.

    (5.5.1) If $h=0$, the diagram of $\Delta$ is
$$\unitlength=1cm
%
$$
For $t=2$, the diagram of
$\tau_{b_2}\tau_{b_1}\tau_{a_0}\tau_{a_1}\tau_c\tau_{a_0}\tau_{a_t}\cdots\tau_{a_1}
\tau_{b_1}\tau_{a_0}\tau_{b_2}\tau_{b_1}\tau_c\tau_{a_0}\tau_{a_1}\cdots\tau_{a_t}\Delta$
is
$$\unitlength=1cm
%
$$
For $t=3$, the diagram of
$\tau_{a_3}\tau_{a_2}\tau_{a_1}\tau_{a_0}\tau_{b_1}\tau_c\tau_{a_0}\tau_{a_1}\tau_{b_2}\tau_{b_1}
\tau_{a_0}\tau_{a_2}\tau_{a_1}\tau_c\tau_{a_0}\tau_{a_t}\cdots\tau_{a_1}\tau_{b_1}\\
\tau_{a_0}\tau_{b_2}\tau_{b_1}\tau_c\tau_{a_0}\tau_{a_1}\cdots
\tau_{a_t}\Delta$ is
$$\unitlength=1cm
%
$$
For $t=4$, the diagram of $\tau_{a_4}\tau_{a_3}\tau_{a_2}\tau_{a_1}
\tau_{a_0}\tau_{b_1}\tau_c\tau_{a_0}\tau_{a_1}\tau_{b_2}
\tau_{b_1}\tau_{a_0}
\tau_{a_2}\tau_{a_1}\tau_c\tau_{a_0}\tau_{a_3}\tau_{a_2}\tau_{a_1}\\
\tau_{b_1} \tau_{a_0}\tau_{a_4}\tau_{a_3}\tau_{a_2}\tau_{a_1}
\tau_c\tau_{b_2}\tau_{b_1}
\tau_{a_0}\tau_{a_1}\tau_{b_1}\tau_c\tau_{a_0}\tau_{a_2}\tau_{a_1}\tau_{b_2}
\tau_{b_1}\tau_{a_0}\tau_{a_3}\tau_{a_2}\tau_{a_1}\tau_c\tau_{a_0}\tau_{a_t}
\cdots\tau_{a_1}\tau_{b_1}\tau_{a_0}\\
\tau_{b_2}\tau_{b_1}\tau_c\tau_{a_0}\tau_{a_1}\cdots\tau_{a_t}\Delta$
is
$$\unitlength=1cm
$$

    (5.5.2) If $h=1$, then we have four cases. The first one is
$$\unitlength=1cm
%
$$
The diagram of $\tau_{a_t}\Delta$ is
$$\unitlength=1cm
%
$$

The second one is
$$\unitlength=1cm
%
$$
The diagram of
$\tau_{b_2}\tau_{b_1}\tau_{a_0}\tau_{a_1}\tau_{a_2}\cdots\tau_{a_t}\Delta$
is
$$\unitlength=1cm
%
$$

The third one is
$$\unitlength=1cm
%
$$
The diagram of $\tau_c\tau_{a_0}\tau_{a_1}\tau_{a_2}\Delta$ is
$$\unitlength=1cm
%
$$

The last one is
$$\unitlength=1cm
%
$$
The diagram of
$\tau_{a_t}\cdots\tau_{a_1}\tau_{a_0}\tau_{b_1}\tau_c\tau_{a_0}\tau_{a_1}\tau_{a_2}
\cdots\tau_{a_t}\Delta$ is
$$\unitlength=1cm
%
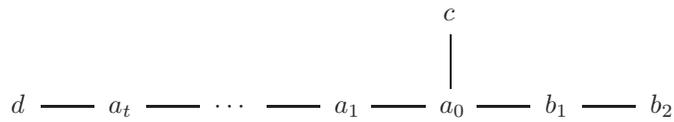
$$

    (5.6) Suppose $D(A)\setminus\{d\}$ is of type $F_4$. Also we have the following cases.

    (5.6.1) If $h=0$, the diagram of $\Delta$ is
$$\unitlength=1cm
%
$$
The diagram of $\tau_{a_1}\tau_{a_0}\tau_{b_1}\tau_{b_2}\tau_{a_0}
\tau_{b_1}\tau_{a_1}\tau_{a_0}\tau_{b_1}\tau_{b_2}
\tau_{a_1}\tau_{a_0}\tau_{b_1}\tau_{a_0}\tau_{a_1}\Delta$ is
$$\unitlength=1cm
%
$$

    (5.6.2) If $h=1$, then there are two cases. One is
$$\unitlength=1cm
%
$$
The diagram of $\tau_{a_1}\Delta$ is
$$\unitlength=1cm
%
$$
Another one is
$$\unitlength=1cm
%
$$
The diagram of
$\tau_{a_1}\tau_{a_0}\tau_{b_1}\tau_{a_0}\tau_{a_1}\Delta$ is
$$\unitlength=1cm
%
$$

    (5.7) Finally, we assume that $D(A)\setminus\{d\}$ is of type $G_2$.

    (5.7.1) If $h=0$, the diagram of $\Delta$ is
$$\unitlength=1cm
%
$$
The diagram of $\tau_{a_0}\tau_b\tau_{a_0}\tau_b\tau_{a_0}\Delta$ is
$$\unitlength=1cm
%
$$

    (5.7.2) If $h=1$, the diagram of $\Delta$ is only one
    possibility:
$$\unitlength=1cm
%
$$
The diagram of $\tau_{a_0}\Delta$ is
$$\unitlength=1cm
%
$$

    {\it Case 6}. Assume that $A_{da_t}=1$ and $A_{a_td}=2$.

    (6.1) Suppose $D(A)\setminus\{d\}$ is of type  $A_l$.

    (6.1.1) When $h=0$, the diagram of $\Delta$ is
$$\unitlength=1cm
%
$$
The diagram of $\tau_{a_1}\tau_{a_2}\cdots\tau_{a_t}\Delta$ is
$$\unitlength=1cm
%
$$

    (6.1.2) When $h=1$, there are two cases. One is
$$\unitlength=1cm
%
$$
The diagram of $\tau_{a_{t}}\Delta$ is
$$\unitlength=1cm
%
$$
Another one is
$$\unitlength=1cm
%
$$
The diagram of $\tau_{a_{2}}\tau_{a_{3}}\Delta$ is
$$\unitlength=1cm
%
$$

    (6.2) Suppose $ D(A)\setminus\{d\}$ is of type $B_l$.

    (6.2.1) When $h=0$, the diagram of $\Delta$ is
$$\unitlength=1cm
%
$$
The diagram of
$\tau_{a_{t}}\cdots\tau_{a_0}\tau_b\tau_{a_{0}}\cdots\tau_{a_t}\Delta$
is
$$\unitlength=1cm
%
$$

    (6.2.2) When $h=1$, there are two cases. One is
$$\unitlength=1cm
%
$$
The diagram of $\tau_{a_{t}}\Delta$ is
$$\unitlength=1cm
%
$$
Another one is
$$\unitlength=1cm
%
$$
The diagram of $\tau_{a_{0}}\Delta$ is
$$\unitlength=1cm
%
$$

    (6.3) Suppose $D(A)\setminus\{d\}$ is of type $C_l$.

    (6.3.1) When $h=0$, the diagram of $\Delta$ is
$$\unitlength=1cm
%
$$
The diagram of
$\tau_{a_{t}}\cdots\tau_{a_0}\tau_b\tau_{a_{0}}\cdots\tau_{a_t}\Delta$
is
$$\unitlength=1cm
%
$$

    (6.3.2) When $h=1$, there are two cases. One is
$$\unitlength=1cm
%
$$
The diagram of $\tau_{a_{t}}\Delta$ is
$$\unitlength=1cm
%
$$
The diagram of $\tau_b\tau_{a_0}\cdots\tau_{a_t}\Delta$ is
$$\unitlength=1cm
%
$$

    (6.3.3) When $h=2$, the diagram of $\Delta$ is only one
    possibility:
$$\unitlength=1cm
%
$$
The diagram of $\tau_{a_{0}}\Delta$ is
$$\unitlength=1cm
%
$$

    (6.4) Suppose $D(A)\setminus\{d\}$ is of type $D_l$.

    (6.4.1) When $h=0$, the diagram of $\Delta$ is
$$\unitlength=1cm
%
$$
The diagram of
$\tau_{a_t}\cdots\tau_{a_0}\tau_c\tau_b\tau_{a_0}\cdots\tau_{a_t}\Delta$
is
$$\unitlength=1cm
%
$$
The diagram of $\tau_{a_t}\Delta$ is
$$\unitlength=1cm
%
$$
The diagram of $\tau_b\tau_{a_0}\tau_{a_1}\Delta$ is
$$\unitlength=1cm
%
$$

    (6.5) Finally, $D(A)\setminus\{d\}$ can not be of type $E_l, F_4,\mbox{ or }  G_2$.\\

    {\it Case 7}. Assume that $A_{da_t}=1$ and  $A_{a_td}=3$.

    (7.1) Suppose $D(A)\setminus\{d\}$ is of type $A_l$.

    It is easy to see that $D(A)\setminus\{d\}$ must be of type $A_1$ or $A_2$.

    If it is $A_1$, $\tau_{a_1}\Delta$ is a GCM-root basis of type
$G_2$.

    If it is $A_2$, there are only two cases. One is
$$\unitlength=1cm
%
$$
The diagram of $\tau_{a_1}\tau_{a_2}\Delta$ is
$$\unitlength=1cm
%
$$
The diagram of $\tau_{a_2}\Delta$ is
$$\unitlength=1cm
%
$$

    (7.2) $D(A)\setminus\{d\}$ can not be of type $B_l, C_l, D_l, E_6, E_7, E_8, F_4$, or $G_2$.

    {\it Case 8}. Assume that $A_{da_t}=1$ and $A_{a_td}=4$.

    (8.1) Suppose $D(A)\setminus\{d\}$ is of type $A_l$.

    It is clear that $D(A)\setminus\{d\}$ must be of type $A_1$.
Thus $\tau_{a_1}\Delta$ is a GCM-root basis of type $A^{(2)}_2$.

    (8.2) Suppose $D(A)\setminus\{d\}$ is $C_l$. We can easily see that there is only one case:\\
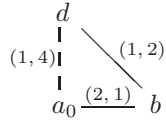
$$\unitlength=1cm
%
$$
The diagram of $\tau_b\tau_{a_0}\Delta$ is
$$\unitlength=1cm
%
$$

    (8.3) $D(A)\setminus\{d\}$ can not be of type $B_l(l>2), D_l, E_6, E_7, E_8, F_4$, or $G_2$.

     Up to now, we have discussed all possible cases and the proof  is
     finished.
\end{proof}

Now we can show  one of our main results.

\begin{theorem}\label{SIM-equiv-GCM}
Let $A=(A_{ij})_{i,j\in I}$ be an SIM such that $A$ is positive
semi-definite and $\corank A\leq 1$. Then for any $d\in I$ such that
$(A_{ij})_{i,j\in I\setminus\{d\}}$ is positive definite, one can
use APR-transformations at some vertices in $I\setminus\{d\}$ such
that $A$ is braid-equivalent to a GCM.
\end{theorem}

\begin{proof}
By induction on $|I|$, we shall prove that for any $d\in I$, if
$(A_{ij})_{i,j\in I\setminus\{d\}}$ is positive definite, then
$D(A)$ is braid-equivalent to a solid diagram by  successive
APR-transformations at some vertices in $I\setminus \{d\}$.

    Obviously $|I|>1$, and there exists  $d\in I$
such that $(A_{ij})_{i,j\in I\setminus \{d\}}$ is positive definite.
Then by the induction hypothesis, the SIM-root basis $\Delta$ of $A$
is braid-equivalent to an SIM-root basis $\Delta'$ with structural
matrix $A'=(A'_{ij})_{i, j\in I}$ such that $D(A')\setminus\{d\}$ is
a disjoint union of some classic Dynkin diagrams by successive
APR-transformations at some vertices in $I\setminus \{d\}$.

    If $D(A')\setminus\{d\}$ is not connected, then there exist two
nonempty and disjoint subsets $I_1$ and $I_2$ of $I\setminus\{d\}$
such that $I\setminus\{d\}=I_1\dot{\cup}I_2$ and
$D(A')\setminus\{d\} =D(A_1')\dot{\cup}D(A_2')$, where
$A_t'=(A_{ij}')_{i,j\in I_t}$, $t=1, 2$. Thus there exist
$a_1,\cdots, a_m\in I_1$ and $b_1,\cdots, b_n\in I_2$ such that both
$\tau_{a_m}\cdots\tau_{a_1}(D(A'), I_1\cup\{d\})$ and
$\tau_{b_n}\cdots\tau_{b_1}(D(A'), I_2\cup\{d\})$  are solid
diagrams, where $(D(A'), I_t\cup \{d\})$ is the full subdiagram of
$D(A')$ with vertex sets $I_t\cup\{d\}$, for $t=1,2$. Note that
$D(A')\setminus \{d\}$ is a disjoint union of some classic Dynkin
diagrams and so invariant under any APR-transformation at $i\in
I\setminus \{d\}$. Thus $\tau_{b_n}
\cdots\tau_{b_1}\tau_{a_m}\cdots\tau_{a_1}(D(A'))$ is a solid
diagram. Therefore, $D(A)$ is braid-equivalent to a solid diagram by
successive APR-transformations at some vertices in $I\setminus
\{d\}$.

    If $D(A')\setminus\{d\}$ is  connected, then $D(A')\setminus
\{d\}$ is the Dynkin diagram of type $A_l, B_l, C_l, D_l, E_6, E_7,
E_8, F_4, G_2$. Then by the induction hypothesis,  $D(A')$ is
braid-equivalent to $D(A'')$ by the APR-transformations at some
vertices in $I\setminus\{d, a_t\}$, where $A''=(A_{ij}'')_{i,j\in
I}$ is an SIM such that $D(A'')\setminus\{a_t\}$ is  a solid
diagram. It is clear that  $D(A'')\setminus\{d\}=D
(A')\setminus\{d\}$. Therefore, by Lemma~\ref{SIM-GCM}, $D(A'')$ is
braid-equivalent to a solid diagram by successive
APR-transformations at some vertices in $I\setminus\{d\}$. This
implies that $D(A)$ is also braid-equivalent to a solid diagram by
successive APR-transformations at some vertices in
$I\setminus\{d\}$.
\end{proof}


Let us recall the notion and some results on IM-Lie algebras of
Slodowy in  \cite{Sl1} and \cite{Sl2}.

Let $A=(A_{ij})_{i,j\in I}$ be an IM, $\Delta$ be the IM-root basis
with structure matrix $A$, and $(-,\,-)$ be the symmetric bilinear
form induced by $A$.

We first introduce an `auxiliary Lie algebra' $\widetilde{{\mathfrak
g}}(A)$, whose generators are $e_\alpha, \alpha\in\pm\Delta$ and
${\mathfrak h}=H\otimes_{\mathbb{Q}}{\mathbb{C}}$, and whose
defining relations are as follows:

\begin{itemize}
\item[(A1)] \,$[h, h']=0$, if $ h, h' \in {\mathfrak h}$;

\item[(A2)] \,$[h, e_\alpha]=(h,\,\alpha)e_\alpha$, if $ h\in
{\mathfrak h}, \alpha \in \pm \Delta$;

\item[(A3)] \,$[e_{\alpha}, e_{-\alpha}]=-\alpha^\vee$, if
$\alpha\in\pm\Delta$.
\end{itemize}

Thus, by defining $\deg e_\alpha=\alpha$ and $\deg h=0$, the root
lattice $\Gamma={\mathbf Z}\cdot\Delta$ produces a gradation:
$\widetilde{{\mathfrak g}}(A)=\bigoplus\limits_{\gamma\in
\Gamma}{\widetilde{\mathfrak g}}_\gamma$.

\begin{definition}\label{IM-Lie-alg}
{\rm  (See  \cite{Sl1}, \cite{Sl2}) \rm Let  $A$  be an intersection
matrix, and ${\mathfrak r}$ be the ideal of $\widetilde{{\mathfrak
g}}(A)$ generated by the following subspaces:
\begin{itemize}
\item[(IM4)] \ $\widetilde{\mathfrak g}_\gamma$, where $\gamma\in\Gamma$ and  $(\gamma, \gamma)>2$.
\end{itemize}

We call the quotient Lie algebra IM$(A)=\widetilde{{\mathfrak
g}}(A)/{\mathfrak r}$ the {\it intersection matrix Lie algebra}
associated to $A$, or  IM-Lie algebra for brevity.}
\end{definition}

The following two results belong to Slodowy.

\begin{proposition}\label{IM(A)-IM(A')}
Let $\Delta$ and $\Delta'$ be two IM-bases with structure matrices
$A$ and $A'$ respectively. If they are braid-equivalent, then
IM$(A)\simeq$IM$(A')$. \hfill$\square$
\end{proposition}

\begin{proposition}\label{SIM=g(A)} Let  $A$ be a
symmetric generalized Cartan matrix. Then IM$(A)$ is isomorphic to
the Kac-Moody Lie algebra ${\mathfrak g}(A)$. \hfill$\square$
\end{proposition}

By Proposition~\ref{IM(A)-IM(A')}, Proposition~\ref{SIM=g(A)} and
Theorem~\ref{SIM-equiv-GCM}, we have the following consequence.

\begin{corollary}\label{IM-KM} Let $A$ be a positive semi-definite IM
with $\corank A\leq 1$. Then IM$(A)$ is isomorphic to a Kac-Moody
Lie algebra. \hfill$\square$
\end{corollary}

\section{$d$-fold Affinizations and Classification of Positive Semi-Definite SIM's}

    In this section, we classify all positive semi-definite SIM's up to braid-equivalence. For
this aim, we first introduce the notion of $d$-fold affinizations.

\subsection{$d$-fold Affinizations}

Let $C$ be an $l\times l$ Cartan matrix, $\dot R$ be the root system
of $C$, and
 $\dot\Delta={\{\alpha_1, \alpha_2,
\cdots, \alpha_l}\}\subseteq \dot R$ be the root basis. Let
$(-,\,-)_C$ be the symmetric bilinear form induced by $C$ (so we
have $C_{ij}=\frac{2(\alpha_i,\;\alpha_j)_C}
{(\alpha_i,\;\alpha_i)_C} $). We define the $d$-fold affinizations
of $C$ as follows:
\begin{definition}\label{affinization}
{\rm   {\it A $d$-fold affinization of $C$} is an $(l+d)$-matrix
$C^{[d]}$ of the form
$$C^{[d]}_{ij}=\frac{2(\alpha_i,\;\alpha_j)_C}{(\alpha_i,\;\alpha_i)_C},
\quad \mbox{where}\quad\alpha_{l+1}, \cdots, \alpha_{l+d}\in\dot
R\cup 2\dot R,$$ such that all $C^{[d]}_{ij}$ are integers. Here
$2\dot R=\{2\alpha\, |\, \alpha\in \dot R\}$.

Particularly, if the above $\alpha_{l+1}, \cdots,
\alpha_{l+d}\in\dot\Delta\cup2\dot\Delta$, then we call the
corresponding $d$-fold affinization matrix $C^{[d]}$ {\it
standard}.}
\end{definition}

\begin{remark}
{\rm Our definition of $d$-fold affinization matrices is a
generalization of that given by S. Berman and R. V. Moody in
\cite{BM} and by G. Benkart and E. Zelmanov in \cite{BZ}. They only
considered the case that $\alpha_{l+1}, \cdots, \alpha_{l+d}\in\dot
R$.}
\end{remark}

\begin{proposition}\label{aff-symm-pos-semi-def}
Let $C$ be a Cartan matrix and $C^{[d]}$ is a $d$-fold affinization
of $C$. Then $C^{[d]}$ is a symmetrizable positive semi-definite
matrix with $\corank C^{[d]}=d$.
 \end{proposition}

\begin{proof}
It is a direct checking.
\end{proof}

In Definition~\ref{affinization}, for the $d$-fold affinization
matrix $C^{[d]}$, we have the SIM-root basis $\Delta=\{\alpha_1,
\alpha_2, \cdots, \alpha_l,\beta_{l+1}, \cdots, \beta_{l+d} \}$,
where $\beta_{l+1}, \cdots, \beta_{l+d}$ correspond respectively to
$\alpha_{l+1}, \cdots, \alpha_{l+d}\in\dot R\cup 2\dot R$, that is
$(\beta_{l+i},\alpha_j)_{C^{[d]}}=(\alpha_{l+i},\alpha_j)_C$,
$(\beta_{l+i},\beta_{l+j})_{C^{[d]}}=(\alpha_{l+i},\alpha_{l+j})_C$
and $(\alpha_{i},\alpha_j)_{C^{[d]}}=(\alpha_{i},\alpha_j)_C$ for
all $i,j$. In this case, we call $\beta_{l+1}, \cdots, \beta_{l+d}$
{\it added roots}; if $\beta_i$ corresponds to $\alpha_i\in \dot R$
(resp. $\alpha_i\in \dot \Delta$), we call
 $\beta_i$ an {\it added single root} (resp. {\it added single simple root});
if $\beta_i$ corresponds to $\alpha_i\in 2\dot R$ (resp.
$\alpha_i\in 2\dot \Delta$), we call $\beta_i$ an {\it added double
root} (resp. {\it added double simple root}).

 Let $C$ be a Cartan matrix,
and let $(-,\,-)_C$ be the induced symmetric bilinear form of $C$.
Then, given a root $\alpha\in\dot R$ we have
$(\alpha,\,\alpha)_C=(\alpha_i,\,\alpha_i)_C$ for some simple root
$\alpha_i$. We call $\alpha$ a {\it short} (resp. {\it long}) {\it
root} if $(\alpha,\,\alpha)=\min_i{(\alpha_i,\,\alpha_i)}$ (resp.
$(\alpha,\,\alpha)=\max_i(\alpha_i,\,\alpha_i)$). Note that if $C$
is symmetric, then all roots have the same square length, which we
also call {\it long roots} by convention. If $C$ is not symmetric,
then every root is either short or long.

Accordingly, $\beta_{l+i}$ is called an added single long (resp.
short) root if $\beta_{l+i}$ corresponds to a long (resp. short)
root $\alpha_{l+i}\in \dot{R}$.

\begin{definition}\label{daff-type}
 {  \rm Let $C$ be an indecomposable Cartan matrix of type $X_l\ (X=A, B, C, D, E, F, G )$, and $\dot
 R$ be the root system of $C$. Let $C^{[d]}$ be a $d$-fold affinization matrix of
 $C$. Let $b$ (resp. $s$) be the number of added single long (resp. short) roots,
 and let $t$ be the number of added double roots.
 Then we call $C^{[d]}$ the $d$-fold affinization matrix of type $X_l^{(b, s, t)}$ ($b+s+t=d$).}
\end{definition}

\begin{lemma}\label{1daf}
Let $C$ be an $l \times l$ Cartan matrix, and $A=(A_{ij})_{i,j\in
I}$ be an $(l+d)$-matrix of rank $l$. Suppose that there is a subset
$J\subset I$ such that $(A_{ij})_{i,j\in J}=C$. Then $A$ is a
$d$-fold affinization matrix $C^{[d]}$ if and only if for all $s\in
I\setminus J$, $(A_{ij})_{i,j\in J \cup \{s\}}$ is a $1$-fold
affinization matrix $C^{[1]}$.
\end{lemma}
\begin{proof} The necessity is obvious. It remains to prove the sufficiency.

For all $s\in I\setminus J$, assume that the added root of the
$1$-fold affinization matrix $(A_{ij})_{i,j\in J \cup \{s\}}$
corresponds to $\alpha_s \in \dot{R}\cup 2\dot{R}$. Let $C^{[d]}$ be
a $d$-fold affinization matrix of $C$,
 such that the added roots $\beta_{l+1},\cdots,\beta_{l+d}$
 correspond respectively to $\alpha_{l+1},\cdots,\alpha_{l+d}$. Then for any $i \in
J$ or $j\in J$ we have $A_{ij}=C^{[d]}_{ij}$, that is, after some
suitable permutations of elements of the index set $I$, if the
matrix $A$ has a block form:
\[
\begin{pmatrix}
C   & A_1 \\
A_2 & A_3
\end{pmatrix}
\]
then the matrix $C^{[d]}$ has a corresponding block form:
\[
\begin{pmatrix}
C   & A_1 \\
A_2 & A'_3
\end{pmatrix}
\]
Since both of them are of rank $l$, we have $A_3=A_2C^{-1}A_1=A'_3$.
So $A=C^{[d]}$.
\end{proof}

\begin{proposition}\label{daff-sdaff}
Let $C$ be a Cartan matrix. Then any $d$-fold affinization matrix
$C^{[d]}$ of $C$ is braid-equivalent to a standard $d$-fold
affinization matrix of the same type.
 \end{proposition}
\begin{proof}
Let the root basis of $C$ be  \ $\dot{\Delta}=\{\alpha_1,\cdots,
\alpha_l\}$, and the SIM-root basis of $C^{[d]}$ be $\Delta
=\{\alpha_1,\cdots, \alpha_l,\beta_{l+1},\cdots,\beta_{l+d}\}$,
where $\beta_{l+1},\cdots,\beta_{l+d}$ are added roots,
corresponding respectively to  $\alpha_{l+1},\cdots,\alpha_{l+d}\in
\dot{R}\cup 2\dot{R}$. Let $\alpha_{l+1}=s_{\alpha_{i_t}}\cdots
s_{\alpha_{i_1}}(\alpha)$, where $\alpha \in \dot{\Delta} \cup
2\dot{\Delta}$, $\alpha_{i_j}\in \dot\Delta$ and $s_{\alpha_{i_j}}$
is a reflection for $j=1,2,\cdots, t$. Denote $\Delta'=
\tau_{i_{t},{l+1}}\Delta$, then $\Delta'=\{ \alpha_1,\cdots,
\alpha_l,s_{\alpha_{i_t}}(\beta_{l+1})
,\beta_{l+2},\cdots,\beta_{l+d}\}$ is an SIM-root basis
braid-equivalent to $\Delta$. Note that
$s_{\alpha_{i_t}}(\beta_{l+1})=
\beta_{l+1}-\frac{2(\alpha_{i_t},\alpha_{l+1})_C}{(\alpha_{i_t},\alpha_{i_t})_C}\alpha_{i_t}$.
It is easy to see that the structural matrix  of $\Delta'$ is also a
$d$-fold affinization and the added roots
$s_{\alpha_{i_t}}(\beta_{l+1})$ correspond to
$s_{\alpha_{i_t}}(\alpha_{l+1})$. By the induction on $t$ we know
that $C^{[d]}$ is braid-equivalent to a $d$-fold affinization
matrix, and its SIM-root basis is $\{\alpha_1,\cdots,
\alpha_l,\gamma_{l+1},\beta_{l+2},\cdots,\beta_{l+d}\}$, such that
the added root $\gamma_{l+1}$ corresponds to $\alpha \in
\dot{\Delta} \cup 2\dot{\Delta}$. Applying the same argument to
$\alpha_{l+2}, \cdots, \alpha_{l+d}$,  we can show that $C^{[d]}$ is
braid-equivalent to a $d$-fold affinization matrix, and its
 SIM-root basis is $\{\alpha_1,\cdots,
\alpha_l,\gamma_{l+1},\gamma_{l+2},\cdots,\gamma_{l+d}\}$, such that
the added roots $\gamma_{l+1},\gamma_{l+2},\cdots,\gamma_{l+d}$
 correspond to the elements of $\dot\Delta\cup
2\dot\Delta$, that is, $C^{[d]}$ is braid-equivalent to a standard
$d$-fold affinization matrix.
\end{proof}

\begin{proposition}\label{aff'form-1'aff'form}
Let $A=(A_{ij})_{i,j\in I}$ be an indecomposable generalized Cartan
matrix of affine type. Then $A$ is a $1$-fold affinization matrix,
and the type of $A$ as a generalized Cartan matrix of affine type
(see the diagrams in \cite{K}) has the following relation with the
type of $A$ as a $1$-fold affinization matrix:\\
{\rm (1)}\ $A$ is of type $D_{4}^{(3)}$ if and only if $A$ is of
type $G_{2}^{(0, 1, 0)}$;\\
{\rm (2)}\ $A$ is of type $G_{2}^{(1)}$ if and only if $A$ is of
type $G_{2}^{(1, 0, 0)}$;\\
{\rm (3)}\ $A$ is of type $E_{6}^{(2)}$ if and only if $A$ is of
type $F_{4}^{(0, 1, 0)}$;\\
{\rm (4)}\ $A$ is of type $F_{4}^{(1)}$ if and only if $A$ is of
type $F_{4}^{(1, 0, 0)}$;\\
{\rm (5)}\ $A$ is of type $D_{l+1}^{(2)}$ if and only if $A$ is of
type $B_{l}^{(0, 1, 0)}$;\\
{\rm (6)}\ $A$ is of type $B_{l}^{(1)}$ if and only if $A$ is of
type $B_{l}^{(1, 0, 0)}$;\\
{\rm (7)}\ $A$ is of type $A_{2l}^{(2)}$ if and only if $A$ is of
type $B_{l}^{(0, 0, 1)}$;\\
{\rm (8)}\ $A$ is of type $A_{2l-1}^{(2)}$ if and only if $A$ is of
type $C_{l}^{(0, 1, 0)}$;\\
{\rm (9)}\ $A$ is of type $C_{l}^{(1)}$ if and only if $A$ is of
type $C_{l}^{(1, 0, 0)}$;\\
{\rm (10)}\ $A$ is of type $A_{l}^{(1)}$ if and only if $A$ is of
type $A_{l}^{(1, 0, 0)}$;\\
{\rm (11)}\ $A$ is of type $A_{2}^{(2)}$ if and only if $A$ is of
type $A_{1}^{(0, 0, 1)}$;\\
{\rm (12)}\ $A$ is of type $D_{l}^{(1)}$ if and only if $A$ is of
type $D_{l}^{(1, 0, 0)}$;\\
{\rm (13)}\ $A$ is of type $E_{l}^{(1)}$ if and only if $A$ is of
type $E_{l}^{(1, 0, 0)}$.
\end{proposition}

\begin{proof}
Suppose $A$ is of type $D_{4}^{(3)}$. Obviously $A$ has an
indecomposable Cartan submatrix $(A_{ij})_{i, j\in J}$ of type
$G_{2}$, where $J$ is a subset of $I$ with $|I|-1$ elements. Denote
 $C=(A_{ij})_{i, j\in J}$ and let $C^{[1]}$ be a $1$-fold affinization matrix of type $G_{2}^{(0, 1, 0)}$.
 By Theorem~\ref{SIM-equiv-GCM}, $C^{[1]}$ is braid-equivalent
to an indecomposable affine generalized Cartan matrix
$A'=(A'_{ij})_{i, j\in I}$ such that the submatrix $(A'_{ij})_{i,
j\in J}=C$. Note that for the $1$-fold affinization matrix
$C^{[1]}$, the added root corresponds to a short root of
 $C$. So the indecomposable affine generalized Cartan matrix $A'$
 has an indecomposable Cartan submatrix $C$ of type $G_{2}$ and another simple root of $A'$
 (which is not a simple root of $C$) is a short
 root. Therefore the type of $A'$ must be $D_{4}^{(3)}$, that is, $A'=A$.
 Hence $A$ is a $1$-fold affinization matrix and the statement (1)
 is true. Similarly for other cases.
\end{proof}

\begin{definition}\label{naff}
{\rm  Let $A'=(A'_{ij})_{i,j\in I}$ be a positive semi-definite SIM,
and $J\subsetneqq I$. We call $A'=(A'_{ij})_{i,j\in I}$ a {\it
$J$-normal affinization matrix}, if $(A'_{ij})_{i,j\in J}$ is a
Cartan matrix, and for any $s\in I\setminus J$, $(A'_{ij})_{i,j\in
J\cup \{ s\}}$ is a positive semi-definite generalized Cartan matrix
of corank $1$.

Let $A=(A_{ij})_{i,j\in I}$ be a positive semi-definite SIM. If $A$
is braid-equivalent to a $J$-normal affinization matrix
$A'=(A'_{ij})_{i,j\in I}$, then we say $A'=(A'_{ij})_{i,j\in I}$
 is a {\it $J$-normal affinization matrix of $A$}.}
\end{definition}

\begin{proposition}\label{ssim-normaff}
Let $A=(A_{ij})_{i,j\in I}$ be a positive semi-definite \SIM\ with
$\rank A= l<|I|$. For any subset $J\subset I$, if $|J|=l$ and $\rank
(A_{ij})_{i,j\in J}=l$, then $A$ is braid-equivalent to a $J$-normal
affinization matrix under a sequence of APR-transformations at the
vertices of $J$.
\end{proposition}

\begin{proof}
Take $s_1\in I\setminus J$, then $(A_{ij})_{i, j\in J\cup\{s_1\}}$
is positive semi-definite of corank $1$. By
Theorem~\ref{SIM-equiv-GCM}, under a sequence of APR-transformations
at the vertices of $J$, $A$ is braid-equivalent to a \SIM\
$A'=(A'_{ij})_{i, j\in I}$, such that $(A'_{ij})_{i, j\in J}$ is a
Cartan matrix and $(A'_{ij})_{i, j\in J\cup\{s_1\}}$ is a positive
semi-definite generalized Cartan matrix of corank $1$. Take $s_2\in
I\setminus (J\cup\{s_1\})$, then similarly under a sequence of
APR-transformations at the vertices of $J$, $A'$ is braid-equivalent
to a \SIM\  $A''=(A''_{ij})_{i, j\in I}$, such that $(A''_{ij})_{i,
j\in J}$ is a Cartan matrix and $(A''_{ij})_{i, j\in J\cup\{s_2\}}$
is a positive semi-definite generalized Cartan matrix of corank $1$.
Note that any symmetrizable generalized Cartan matrix is invariant
under APR-transformations. So we have $(A''_{ij})_{i, j\in
J\cup\{s_1\}}=(A'_{ij})_{i, j\in J\cup\{s_1\}}$. Now the induction
on $|I\setminus J|$ implies the result.
\end{proof}

\begin{lemma}\label{indecom-SIM}
Let $A=(A_{ij})_{i, j\in I}$ be an indecomposable positive
semi-definite \SIM\ of $\rank A=l$. Then there exists an
indecomposable submatrix $(A_{ij})_{i, j\in J}$ such that $|J|=l$
and $\rank (A_{ij})_{i, j\in J}=l$.
\end{lemma}
\begin{proof}
Since $A$ is positive semi-definite and $\rank A=l$, there is a
subset $K\subseteq I$ such that $|K|=l$ andÇÒ$(A_{ij})_{i, j\in K}$
 is a positive definite submatrix. Among all of these subsets $K$, we pick out a subset $J$, such that $(A_{ij})_{i,
j\in J}$ has a minimal number of indecomposable components. We claim
that $(A_{ij})_{i, j\in J}$ is indecomposable.

Otherwise, suppose  $(A_{ij})_{i, j\in J}$ is decomposable. Let
$(A_{ij})_{i, j\in J_1}, \cdots, (A_{ij})_{i, j\in J_m}$ $(m>1)$ be
all indecomposable components of  $(A_{ij})_{i, j\in J}$. By
Proposition~\ref{ssim-normaff}, $A$ is braid-equivalent to a
 $J$-normal affinization matrix $A'=(A'_{ij})_{i, j\in I}$. Obviously, $(A'_{ij})_{i,
j\in J_1}, \cdots, (A'_{ij})_{i, j\in J_m}$ are all indecomposable
components of the Cartan submatrix $(A'_{ij})_{i, j\in J}$, and for
all $s\in I\setminus J$, every indecomposable component of
$(A'_{ij})_{i, j\in J\cup\{s\}}$ is of form $(A'_{ij})_{i, j\in
J_{k_1}\cup\cdots\cup J_{k_n}\cup\{s\}}$, $(A'_{ij})_{i, j\in
J_{k_{n+1}}},\cdots, (A'_{ij})_{i, j\in J_m}$, where $1\leq n\leq
m$, and we also have $(A'_{ij})_{i, j\in J_{k_1}\cup\cdots\cup
J_{k_n}\cup\{s\}}$ is an indecomposable affine generalized Cartan
matrix.

If there is an $s\in I\setminus J$ such that the above $n > 1$, that
is, the number of all indecomposable components of $(A'_{ij})_{i,
j\in J\cup\{s\}}$ is less that $m$, then we let $(A'_{ij})_{i, j\in
J'}$ be the indecomposable Cartan submatrix of the indecomposable
affine generalized Cartan matrix $(A'_{ij})_{i, j\in J_{k_1}\cup
\cdots\cup J_{k_n}\cup\{s\}}$, such that
$|J'|=|J_{k_1}\cup\cdots\cup J_{k_n}\cup\{s\}|-1$. Therefore
$(A'_{ij})_{i, j\in J'\cup J_{k_{n+1}}\cup\cdots\cup J_{k_m}}$ is a
Cartan matrix of order $l$, and the number of its indecomposable
components is $1+(m-n)< m$, which contradicts to the choice of $J$.
Therefore, for any $s\in I\setminus J$, if there is some $1\leq
x\leq m$ such that $(A'_{ij})_{i, j\in J_{k_x}\cup\{s\}}$ is
indecomposable, then $(A'_{ij})_{i, j\in J_{k_x}\cup\{s\}}$ must be
an indecomposable affine generalized Cartan matrix.

Since $A$ is indecomposable, $A'$ is also indecomposable. Then there
exist two submatrices among $(A'_{ij})_{i, j\in J_1}, \cdots,
(A'_{ij})_{i, j\in J_m}$, say  $(A'_{ij})_{i, j\in J_1}$ and
$(A'_{ij})_{i, j\in J_2}$, which are connected with each other by
some vertices in $I\setminus J$, that is, there are
$t_1,\cdots,t_k\in I \ (k\geq 3)$ such that each of $A'_{t_1 t_2},
A'_{t_2 t_3},\cdots, A'_{t_{k-1} t_k}$ is nonzero, where $t_1\in
J_1$, $t_k\in J_2$, but $t_{2},\cdots,t_{k-1}\in I\setminus J$. We
can suppose that such $k$ is minimal. Then for all $a\in J_1$, we
have $A'_{at_{3}}=0$.  By the  result above we know that
$(A'_{ij})_{i, j\in J_{1}\cup\{t_2\}}$ is an indecomposable affine
generalized Cartan matrix. Then it is not difficult to show that
$(A'_{ij})_{i, j\in J_{1}\cup\{t_2, t_{3}\} }$ is (or is
braid-equivalent to) an  indefinite generalized Cartan matrix, as
contradicts to the positive semi-definiteness.
\end{proof}

In the following we come to one of our main results in this section,
which shows that any indecomposable positive semi-definite \SIM\ is
in fact braid-equivalent to a $d$-fold affinization matrix.

\begin{theorem}\label{A-equi-d-fold-aff}
Let $A$ be a positive semi-definite $(l+d)\times (l+d)$ \SIM\ with
$\rank A=l$. Then $A$ is braid-equivalent to a $d$-fold affinization
matrix $C^{[d]}$.
\end{theorem}

\begin{proof}
 It is no harm to assume that $A$ is indecomposable.
 We can choose some $J\subseteq I$ such that $|J|=l$
and $(A_{ij})_{i, j\in J}$ is an indecomposable \SIM\ with $\rank
(A_{ij})_{i, j\in J}=l$, according to Lemma~\ref{indecom-SIM}. So by
Proposition~\ref{ssim-normaff}, $A$ is braid-equivalent to a
$J$-normal affine matrix $(A'_{ij})_{i, j\in I}$. Then the Cartan
submatrix $(A'_{ij})_{i, j\in J}$ is also indecomposable , and for
all $s\in I\setminus J$, $(A'_{ij})_{i, j\in J\cup\{s\}}$ is an
indecomposable affine generalized Cartan submatrix. Now, if
$(A'_{ij})_{i, j\in J}$ is of type $X_l\ (X=A, B, C, D, E, F, G )$,
then we call $J$ a {\it maximal subset} of $A$ of type $X_l$. In the
following we consider all possible cases and choose a suitable
maximal subset $J$, and we prove the corresponding $J$-normal
affinization matrix $(A'_{ij})_{i, j\in I}$ is a $d$-fold
affinization matrix (for the concrete types of affine generalized
Cartan matrices considered here, see the tables of the affine Dynkin
diagrams in \cite{K}).

(1)\, Suppose $A$ has a maximal subset $J$ of type $G_2$. Denote
$(A'_{ij})_{i, j\in I}$ a $J$-normal affinization matrix of $A$.
Then for all $s\in I\setminus J$, the type of the indecomposable
affine generalized Cartan submatrix  $(A'_{ij})_{i, j\in
J\cup\{s\}}$ must be $G_2^{(1)}$ or $D_4^{(3)}$. We denote the
number of $G_2^{(1)}$'s or $D_4^{(3)}$'s respectively by $b$ or $s$.
By Proposition~\ref{aff'form-1'aff'form} and Lemma~\ref{1daf}, we
know $(A'_{ij})_{i, j\in I}$ is a $d$-fold affinization matrix of
type $G_2^{(b, s, 0)}$.

(2)\, Suppose $A$ has a maximal subset $J$ of type $F_4$. Denote
$(A'_{ij})_{i, j\in I}$ a $J$-normal affinization matrix of $A$.
Then for any $s\in I\setminus J$, the type of the indecomposable
affine generalized Cartan submatrix  $(A'_{ij})_{i, j\in
J\cup\{s\}}$ must be $F_4^{(1)}$ or $E_6^{(2)}$. We denote the
number of $F_4^{(1)}$'s or $E_6^{(2)}$'s respectively by $b$ or $s$.
By Proposition~\ref{aff'form-1'aff'form} and Lemma~\ref{1daf}, we
know $(A'_{ij})_{i, j\in I}$ is a $d$-fold affinization matrix of
type $F_4^{(b, s, 0)}$.

(3)\, Suppose $A$ has a maximal subset $J$ of type $B_l\,(l\geq 3)$
and has no   maximal subset of type $F_4$. Denote $(A'_{ij})_{i,
j\in I}$ a $J$-normal affinization matrix of $A$. Then for all $s\in
I\setminus J$, the type of the indecomposable affine generalized
Cartan submatrix $(A'_{ij})_{i, j\in J\cup\{s\}}$ must be
$B_l^{(1)}$,  $D_{l+1}^{(2)}$ or $A_{2l}^{(2)}$. Denote the number
of them  by $b$, $s$ or $t$ respectively. By
Proposition~\ref{aff'form-1'aff'form} and Lemma~\ref{1daf}, we know
that $(A'_{ij})_{i, j\in I}$ is a $d$-fold affinization matrix of
type $B_l^{(b, s, t)}$.

(4)\, Suppose $A$ has a maximal subset $J$ of type $B_2$. Let
$(A'_{ij})_{i, j\in I}$ be a $J$-normal affinization matrix of $A$.
Then for all $s\in I\setminus J$, the type of the indecomposable
affine generalized Cartan submatrix $(A'_{ij})_{i, j\in J\cup\{s\}}$
must be $C_2^{(1)}$, $D_{3}^{(2)}$ or $A_{4}^{(2)}$. Denote the
number of them  by $b$, $s$ or $t$ respectively. By
Proposition~\ref{aff'form-1'aff'form} and Lemma~\ref{1daf}, we know
that $(A'_{ij})_{i, j\in I}$ is a $d$-fold affinization matrix of
type $B_2^{(b, s, t)}$.

(5)\, Suppose $A$ has a maximal subset $J$ of type $C_l\,(l\geq 3)$
and has no  maximal subset of type  $F_4$ or $B_l$. Let
$(A'_{ij})_{i, j\in I}$ be a $J$-normal affinization matrix of $A$.
Then for any $s\in I\setminus J$, the type of the indecomposable
affine generalized Cartan submatrix $(A'_{ij})_{i, j\in J\cup\{s\}}$
must be $C_l^{(1)}$ or $A_{2l-1}^{(2)}$. Denote the number of them
by $b$ or $s$ respectively. By Proposition~\ref{aff'form-1'aff'form}
and Lemma~\ref{1daf}, we know that $(A'_{ij})_{i, j\in I}$ is a
$d$-fold affinization matrix of type $C_l^{(b, s, 0)}$.

(6)\, Suppose $A$ is not symmetric and has no maximal subset of type
$G_2$, $F_4$, $B_l$ or $C_l$. It is easy to see $\rank A=1$. Take
any $j\in I$ such that the length of the root $\alpha_j$ in the
corresponding SIM-root basis is minimal. Denote $J=\{ j\}$ and let
$(A'_{ij})_{i, j\in I}$ be a $J$-normal affinization matrix of $A$.
Then for any $s\in I\setminus J$, the type of the indecomposable
affine generalized Cartan submatrix $(A'_{ij})_{i, j\in J\cup\{s\}}$
must be $A_{1}^{(1)}$ or $A_2^{(2)}$. Denote the number of them by
$b$ or $t$ respectively. By Proposition~\ref{aff'form-1'aff'form}
and Lemma~\ref{1daf}, we know that $(A'_{ij})_{i, j\in I}$ is a
$d$-fold affinization matrix of type $A_1^{(b, 0, t)}$.

We have discussed the case that $A$ is not symmetric as above. In
the following, we consider the case that $A$ is symmetric.

(7)\, Suppose $A$ has a maximal subset $J$ of type $E_l\ (l=6,7,8)$.
 Denote $(A'_{ij})_{i, j\in I}$ a $J$-normal affinization matrix of $A$. Then for any $s\in I\setminus J$,
the type of the indecomposable affine generalized Cartan submatrix
$(A'_{ij})_{i, j\in J\cup\{s\}}$ must be $E_l^{(1)}$. By
Proposition~\ref{aff'form-1'aff'form} and Lemma~\ref{1daf}, we know
$(A'_{ij})_{i, j\in I}$ is a $d$-fold affinization matrix of type
$E_l^{(d, 0, 0)}$.

(8)\, Suppose $A$ has a maximal subset $J$ of type $D_l$ and has no
maximal subset of type $E_l$.
 Denote $(A'_{ij})_{i, j\in I}$ a $J$-normal affinization matrix of $A$. Then for any $s\in I\setminus J$,
the type of the indecomposable affine generalized Cartan submatrix
$(A'_{ij})_{i, j\in J\cup\{s\}}$ must be $D_l^{(1)}$. By
Proposition~\ref{aff'form-1'aff'form} and Lemma~\ref{1daf}, we know
$(A'_{ij})_{i, j\in I}$ is a $d$-fold affinization matrix of type
$D_l^{(d, 0, 0)}$.

(9)\, Suppose $A$ has a maximal subset $J$ of type $A_l$ and has no
maximal subset of other types.
 Denote $(A'_{ij})_{i, j\in I}$ a $J$-normal affinization matrix of $A$. Then for any $s\in I\setminus J$,
the type of the indecomposable affine generalized Cartan submatrix
$(A'_{ij})_{i, j\in J\cup\{s\}}$ must be $A_l^{(1)}$. By
Proposition~\ref{aff'form-1'aff'form} and Lemma~\ref{1daf}, we know
$(A'_{ij})_{i, j\in I}$ is a $d$-fold affinization matrix of type
$A_l^{(d, 0, 0)}$.

We have discussed all possible cases, and in each case, $A$ is
braid-equivalent to a $d$-fold affinization matrix which is at the
same time  a normal affinization matrix.
\end{proof}

\subsection{Classification of Positive Semi-definite Symmetrizable
Intersection Matrices}

As another main result of this section, the following theorem
produces a classification of all indecomposable positive
semi-definite \SIMs by standard $d$-fold affinization matrices, up
to braid-equivalence.

\begin{theorem}\label{classific}
Let $A$ be an indecomposable positive semi-definite \SIM. Then $A$
is braid-equivalent to a standard $d$-fold affinization matrix of
type $X_l^{(b, s, t)}$ ($b+s+t=d$). Moreover, any two standard
$d$-fold affinization matrices with the same type are
braid-equivalent. Here the type $X_l^{(b, s, t)}$ must be one of the
following forms:\\
{\rm (1)}\ $X_l^{(d, 0, 0)}$ for $X_l= A_l(l\geq 2),  \,D_l, \,E_6,
\,E_7, \,E_8$;\\
{\rm (2)}\ $X_l^{(b, s, 0)}$ for $X_l= C_l(l\geq 3), \,F_4, \,G_2$;\\
{\rm (3)}\ $A_1^{(b, 0, t)}$, $B_l^{(b, s, t)}$, where for $B_l^{(b,
s, t)}$, any added double root must correspond to the double of a
short root.
\end{theorem}
\begin{proof}
By Theorem~\ref{A-equi-d-fold-aff} and Proposition~\ref{daff-sdaff},
$A$ is braid-equivalent to a standard $d$-fold affinization matrix.

If $A'=(A'_{ij})_{i, j\in I}$ and $A''=(A''_{ij})_{i, j\in I}$ are
two standard $d$-fold affinization matrices of the same type
$X_l^{(b, s, t)}$, then we can assume that there is a subset
$J\subseteq I$ such that $(A'_{ij})_{i, j\in J}=(A''_{ij})_{i, j\in
J}$ is an indecomposable Cartan matrix of type $X_l$ and for any
$a\in J\setminus I$, $(A'_{ij})_{i, j\in J\cup\{ a\}}$ and
$(A''_{ij})_{i, j\in J\cup\{ a\}}$ are two standard
 $1$-fold affinization matrices of the same type.  So for all $a\in J\setminus
I$, by Proposition~\ref{aff'form-1'aff'form} we see both
$(A'_{ij})_{i, j\in J\cup\{ a\}}$ and $(A''_{ij})_{i, j\in J\cup\{
a\}}$ are braid-equivalent to the same one indecomposable affine
generalized Cartan matrix. Therefore both $A'$ and
 $A''$ are braid-equivalent to the same one $J$-normal affinization matrix. Hence they are
 braid-equivalent to each other.

It is easy to see that there exist standard $d$-fold affinization
matrices of type appearing in (1), (2) and (3). Next we show that
the types appearing in (1), (2) and (3) are the only possibilities.
 For (1) and (2), we only
need to prove the added roots in a standard $d$-fold affinization
matrix $C^{[d]}$ are all single roots, that is, the corresponding
roots $\alpha_{l+1},\cdots,\alpha_{l+d}$ lie in $\dot\Delta$.
Otherwise, suppose some $\alpha_{l+j}=2\alpha_k\in
2\dot\Delta,\,1\leq k\leq l$. Then for all $i=1, \cdots, l$,
$$C^{[d]}_{l+j,\;i}=
\frac{2(\alpha_{l+j},\;\alpha_i)_C}{(\alpha_{l+j},\;\alpha_{l+j})_C}=
\frac{(\alpha_k,\;\alpha_i)_C}{(\alpha_k,\;\alpha_k)_C}=\frac{1}{2}\cdot
C_{ki}.$$ Obviously, for any Cartan matrix $C$ of type $A_l(l\geq
2)$, $C_l(l\geq 3)$, $D_l$, $E_6$, $E_7$, $E_8$, $F_4$ or $G_2$,
there exists an $i$ with $1\leq i\leq l$ such that $C_{ki}=-1$ or
$-3$. Hence $C^{[d]}_{l+1,\;i}=\frac{1}{2}\cdot
C_{ki}\notin\mathbf{Z}$, as contradicts to the definition of
$d$-fold affinizations.

For (3), we only need to prove that for the type $B_l^{(b, s, t)}$,
any added double root must correspond to the double of a short root.
Let $\{\alpha_1,\cdots,\alpha_{l-1},\alpha_l\}$ be a root basis of
$B_l$, where $\alpha_1,\cdots,\alpha_{l-1}$ are the long simple
roots and $\alpha_l$ is the short simple root. Let $C^{[d]}$ be a
standard $d$-affinization matrix of type $B_l^{(b,s,t)}$. Suppose
$\beta_{l+j}$ is an added double root corresponding to
$\alpha_{l+j}=2\alpha_k$ with $1\leq k\leq l-1$, then
$$C^{[d]}_{l+j,\;k+1}=
\frac{2(\alpha_{l+j},\;\alpha_{k+1})_C}{(\alpha_{l+j},\;\alpha_{l+j})_C}=
\frac{(\alpha_k,\;\alpha_{k+1})_C}{(\alpha_k,\;\alpha_k)_C}=\frac{1}{2}\cdot
C_{k,\;k+1}=-\frac{1}{2}.$$ This is impossible and so any added
double root must correspond to the double of a short root.
\end{proof}

\section{Weyl Root Systems of Positive Semi-definite SIM's}

In this section, we give an explicit structure of the Weyl root
system for each  $d$-fold affinization matrices by using the root
system of the Cartan matrix and some special null roots.

Let $A$ be an SIM, and $(-,-)$ be the corresponding symmetric
bilinear form. Set
$${\rm rad}(-,-):=\{\alpha\in H\,|\,
(\alpha,\,\beta)=0\quad\mbox{for all}\quad\beta\in
  H\}$$ and
  $${\rm rad}(\Gamma):={\rm rad}(-,-)\cap \Gamma,$$
where $\Gamma$ is the root lattice. We call each element in ${\rm
rad}(\Gamma)$ a {\it null root}.

For the sake of convenience, we relabel the simple roots in Dynkin
diagrams as follows:

$A_l$:
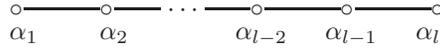
$$\unitlength=1cm
\begin{picture}(9,1)
\put(0,0.3){$\circ$} \put(0,0){$\alpha_1$}
\put(0.2,0.4){\line(1,0){1}} \put(1.2,0.3){$\circ$}
\put(1.2,0){$\alpha_2$} \put(1.4,0.4){\line(1,0){0.6}}
\put(2.1,0.3){$\cdots$} \put(2.6,0.4){\line(1,0){0.6}}
\put(3.2,0.3){$\circ$} \put(3,0){$\alpha_{l-2}$}
\put(3.4,0.4){\line(1,0){1}} \put(4.4,0.3){$\circ$}
\put(4.2,0){$\alpha_{l-1}$} \put(4.6,0.4){\line(1,0){1}}
\put(5.6,0.3){$\circ$} \put(5.4,0){$\alpha_{l}$}
\end{picture}$$

$B_l\;(l\geq 3)$:
$$\unitlength=1cm
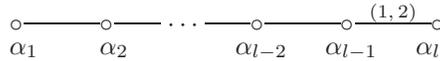
\begin{picture}(9,1)
\put(0,0.3){$\circ$} \put(0,0){$\alpha_1$}
\put(0.2,0.4){\line(1,0){1}} \put(1.2,0.3){$\circ$}
\put(1.2,0){$\alpha_2$} \put(1.4,0.4){\line(1,0){0.6}}
\put(2.1,0.3){$\cdots$} \put(2.6,0.4){\line(1,0){0.6}}
\put(3.2,0.3){$\circ$} \put(3,0){$\alpha_{l-2}$}
\put(3.4,0.4){\line(1,0){1}} \put(4.4,0.3){$\circ$}
\put(4.2,0){$\alpha_{l-1}$}
\put(5.1,0.55){\makebox(0,0)[c]{\scriptsize $(1,2)$}}
\put(4.6,0.4){\line(1,0){1}} \put(5.6,0.3){$\circ$}
\put(5.4,0){$\alpha_{l}$}
\end{picture}$$

$C_l$:
$$\unitlength=1cm
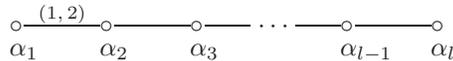
\begin{picture}(9,1)
\put(0,0){$\alpha_1$} \put(0,0.3){$\circ$}
\put(0.2,0.4){\line(1,0){1}}
\put(0.7,0.55){\makebox(0,0)[c]{\scriptsize $(1,2)$}}
\put(1.2,0.3){$\circ$} \put(1.2,0){$\alpha_2$}
\put(1.4,0.4){\line(1,0){1}} \put(2.4,0.3){$\circ$}
\put(2.4,0){$\alpha_3$} \put(2.6,0.4){\line(1,0){0.6}}
\put(3.3,0.3){$\cdots$} \put(3.8,0.4){\line(1,0){0.6}}
\put(4.4,0.3){$\circ$} \put(4.4,0){$\alpha_{l-1}$}
\put(4.6,0.4){\line(1,0){1}} \put(5.6,0.3){$\circ$}
\put(5.6,0){$\alpha_l$}
\end{picture}$$

$D_l$:
$$\unitlength=1cm
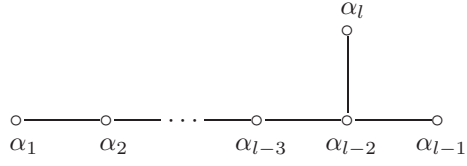
\begin{picture}(9,2)
\put(0,0.3){$\circ$} \put(0,0){$\alpha_1$}
\put(0.2,0.4){\line(1,0){1}} \put(1.2,0.3){$\circ$}
\put(1.2,0){$\alpha_2$} \put(1.4,0.4){\line(1,0){0.6}}
\put(2.1,0.3){$\cdots$} \put(2.6,0.4){\line(1,0){0.6}}
\put(3.2,0.3){$\circ$} \put(3,0){$\alpha_{l-3}$}
\put(3.4,0.4){\line(1,0){1}} \put(4.4,0.3){$\circ$}
\put(4.2,0){$\alpha_{l-2}$} \put(4.6,0.4){\line(1,0){1}}
\put(5.6,0.3){$\circ$} \put(5.4,0){$\alpha_{l-1}$}
\put(4.5,0.5){\line(0,1){1}} \put(4.4,1.5){$\circ$}
\put(4.4,1.8){$\alpha_l$}
\end{picture}$$

$E_6$:
$$\unitlength=1cm
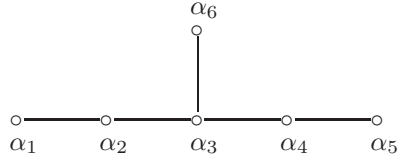
\begin{picture}(9,2)
\put(0,0){$\alpha_1$} \put(0,0.3){$\circ$}
\put(0.2,0.4){\line(1,0){1}} \put(1.2,0.3){$\circ$}
\put(1.2,0){$\alpha_2$} \put(1.4,0.4){\line(1,0){1}}
\put(2.4,0.3){$\circ$} \put(2.4,0){$\alpha_3$}
\put(2.6,0.4){\line(1,0){1}} \put(3.6,0.3){$\circ$}
\put(3.6,0){$\alpha_4$} \put(3.8,0.4){\line(1,0){1}}
\put(4.8,0.3){$\circ$} \put(4.8,0){$\alpha_5$}
\put(2.5,0.5){\line(0,1){1}} \put(2.4,1.5){$\circ$}
\put(2.4,1.8){$\alpha_6$}
\end{picture}$$

$E_7$:
$$\unitlength=1cm
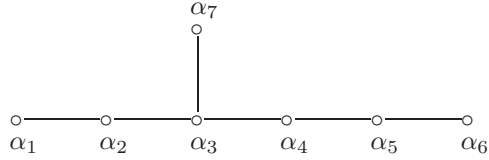
\begin{picture}(9,2)
\put(0,0){$\alpha_1$} \put(0,0.3){$\circ$}
\put(0.2,0.4){\line(1,0){1}} \put(1.2,0.3){$\circ$}
\put(1.2,0){$\alpha_2$} \put(1.4,0.4){\line(1,0){1}}
\put(2.4,0.3){$\circ$} \put(2.4,0){$\alpha_3$}
\put(2.6,0.4){\line(1,0){1}} \put(3.6,0.3){$\circ$}
\put(3.6,0){$\alpha_4$} \put(3.8,0.4){\line(1,0){1}}
\put(4.8,0.3){$\circ$} \put(4.8,0){$\alpha_5$}
\put(5,0.4){\line(1,0){1}} \put(6,0.3){$\circ$}
\put(6,0){$\alpha_6$} \put(2.5,0.5){\line(0,1){1}}
\put(2.4,1.5){$\circ$} \put(2.4,1.8){$\alpha_7$}
\end{picture}$$

$E_8$:
$$\unitlength=1cm
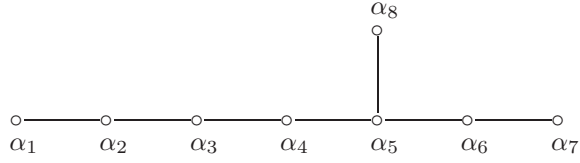
\begin{picture}(9,2)
\put(0,0){$\alpha_1$} \put(0,0.3){$\circ$}
\put(0.2,0.4){\line(1,0){1}} \put(1.2,0.3){$\circ$}
\put(1.2,0){$\alpha_2$} \put(1.4,0.4){\line(1,0){1}}
\put(2.4,0.3){$\circ$} \put(2.4,0){$\alpha_3$}
\put(2.6,0.4){\line(1,0){1}} \put(3.6,0.3){$\circ$}
\put(3.6,0){$\alpha_4$} \put(3.8,0.4){\line(1,0){1}}
\put(4.8,0.3){$\circ$} \put(4.8,0){$\alpha_5$}
\put(5,0.4){\line(1,0){1}} \put(6,0.3){$\circ$}
\put(6,0){$\alpha_6$} \put(6.2,0.4){\line(1,0){1}}
\put(7.2,0.3){$\circ$} \put(7.2,0){$\alpha_7$}
\put(4.9,0.5){\line(0,1){1}} \put(4.8,1.5){$\circ$}
\put(4.8,1.8){$\alpha_8$}
\end{picture}$$

$F_4$:
$$\unitlength=1cm
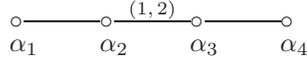
\begin{picture}(9,1)
\put(0,0.3){$\circ$} \put(0,0){$\alpha_1$}
\put(0.2,0.4){\line(1,0){1}} \put(1.2,0.3){$\circ$}
\put(1.2,0){$\alpha_2$} \put(1.9,0.55){\makebox(0,0)[c]{\scriptsize
$(1,2)$}} \put(1.4,0.4){\line(1,0){1}} \put(2.4,0.3){$\circ$}
\put(2.4,0){$\alpha_3$} \put(2.6,0.4){\line(1,0){1}}
\put(3.6,0.3){$\circ$} \put(3.6,0){$\alpha_{4}$}
\end{picture}$$

$G_2$:
$$\unitlength=1cm
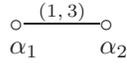
\begin{picture}(9,1)
\put(0,0.3){$\circ$} \put(0,0){$\alpha_1$}
\put(0.2,0.4){\line(1,0){1}} \put(1.2,0.3){$\circ$}
\put(1.2,0){$\alpha_2$} \put(0.7,0.55){\makebox(0,0)[c]{\scriptsize
$(1,3)$}}
\end{picture}$$

In the above labeling , $\alpha_{1}$ is always a long root and
$\alpha_{l}$(where $l$ is the maximal number in the labeling) is
always a short root.

Let $A$ be an SIM which is positive semi-definite. To investigate
its Weyl root system $R^W$, we can assume that $A$ is a standard
$d$-fold affinization matrix of type $X_l^{(b, s, t)}$. We can
further assume that $A$ is so `standard' that each added long root
corresponds to $\alpha_{1}$, each added short root corresponds to
$\alpha_{l}$, and each added double short root corresponds to
$2\alpha_{l}$.

Therefore we can assume that all added long roots corresponding to
$\alpha_1$ are $\alpha_{1j}, j=1, 2, \cdots, b$; all added short
roots corresponding to $\alpha_l$ are $\alpha_{lj}, j=1, 2, \cdots,
s$; and all added double short roots corresponding to $2\alpha_l$
are $\beta_{lj}, j=1, 2, \cdots, t$. Denote
$\delta_{1j}:=\alpha_{1j}-\alpha_1,\,j=1,\cdots,b$,
$\delta_{lj}:=\alpha_{lj}-\alpha_l,\,j=1,\cdots,s$, and
$\delta'_{lj}:=\beta_{lj}-2\alpha_l,\,j=1,\cdots,t$. It is easy to
see that $\{\delta_{11}, \cdots, \delta_{1b}, \delta_{l1}, \cdots,
\delta_{ls}, \delta'_{l1}, \cdots, \delta'_{lt}\}$ is a basis of
rad$(-, -)$.

By the classification theorem (Theorem~\ref{classific}) we know that
the  form of the type $X_l^{(b, s, t)}$ is: (1) $X_l^{(d, 0, 0)}$,
where $X_l= A_l(l\geq 2), \,D_l, \,E_6, \,E_7, \,E_8$; (2) $X_l^{(b,
s, 0)}$, where $X_l= C_l(l\geq 3), \,F_4, \,G_2$; or (3) $A_1^{(b,
0, t)}$, $B_l^{(b, s, t)}$. Thus the following theorem gives a
complete description of the Weyl root system $R^W$ of $A$ in terms
of the root system $\dot{R}$ of the Cartan matrix of type $X_l$ and
the null roots $\delta_{11}, \cdots, \delta_{1b}, \delta_{l1},
\cdots, \delta_{ls}, \delta'_{l1}, \cdots, \delta'_{lt}$. Note that
these Weyl roots have at most three lengths. So we denote all short
roots, all middle roots and all long roots respectively by $R_s^W$,
$R_m^W$ and $R_l^W$.

\begin{theorem}\label{weyl-root} The assumptions are as above. Then

${\mathrm (1)}$\ For the type $X_l^{(d, 0, 0)}$, where $X_l=
A_l(l\geq 2), \,D_l, \,E_6, \,E_7, \,E_8$, we have
$$R^{W}=\dot{R}+rad(\Gamma).$$

$\mathrm{(2)}$\ For the type $X_l^{(b, s, 0)}$, where $X_l=
C_l(l\geq 3), \,F_4, \,G_2$, we have $$R^W=R_s^W\bigcup R_l^W,$$
where
$$R_s^W=\dot{R}_s+rad(\Gamma)$$ and $R_l^W$ is the following set:

$\mathrm{(2.1)}$\ for $X_l=C_l (l\geq 3)$,
\begin{eqnarray*}
R^{W}_l & = & \{\alpha+\sum\limits_{j=1}^b
a_{j}\delta_{1j}+2\sum\limits_{j=1}^s b_{j}\delta_{lj}\;|\;
\alpha\in\dot{R}_l, a_{j}, b_j\in {\mathbf Z},\\
 & & \mbox{ and there is
at most one odd }a_j,j=1,\cdots,b\};
\end{eqnarray*}

$\mathrm{(2.2)}$\ for $X_l=F_4$,
$$R^{W}_l=\{\alpha+\sum\limits_{j=1}^b a_{j}\delta_{1j}+2\sum\limits_{j=1}^s b_{j}\delta_{lj}\;|\;
\alpha\in\dot{R}_l, a_{j}, b_j\in {\mathbf Z}\};$$

$\mathrm{(2.3)}$\ for $X_l=G_2$,
$$R^{W}_l=\{\alpha+\sum\limits_{j=1}^b a_{j}\delta_{1j}+3\sum\limits_{j=1}^s b_{j}\delta_{lj}\;|\;
\alpha\in\dot{R}_l, a_{j}, b_j\in {\mathbf Z}\}.$$

$\mathrm{(3)}$\ For the type $A_1^{(b, 0, t)}$, we have
$R^W=R_s^W\bigcup R_l^W$, where
\begin{eqnarray*} R^{W}_s & = &
\{\pm\alpha_1+\sum\limits_{j=1}^b
a_{j}\delta_{1j}+\sum\limits_{j=1}^t b_{j}
\delta'_{1j}\;|\;a_{j},b_j\in {\mathbf Z}, \\
 & & \mbox{ and there is at
most one odd }a_j,j=1\cdots,b \}.
\end{eqnarray*}
and
\begin{eqnarray*}
R^{W}_l & = & \{\pm2\alpha_1+4\sum\limits_{j=1}^b
a_{j}\delta_{1j}+\sum\limits_{j=1}^t b_{j}\delta'_{1j}\;|\;
a_{j},b_j\in {\mathbf Z},\\
 & &  \mbox{ and there is exactly one odd }
b_j,j=1,\cdots,t \}.
\end{eqnarray*}

$\mathrm{(4)}$\ For the type $B_l^{(b, s, t)}$, we have
$R^W=R_s^W\bigcup R_m^W\bigcup R_l^W$, where
\begin{eqnarray*}
R^{W}_s & = & \{\alpha+\sum\limits_{j=1}^b
a_{j}\delta_{1j}+\sum\limits_{j=1}^s b_{j}\delta_{lj}
+\sum\limits_{j=1}^t c_{j}\delta'_{lj}\;|\;\alpha\in \dot{R}_s,
a_{j},b_j, c_j\in {\mathbf Z},
 \\
 & &  \mbox{and there is at most one odd} \ b_j,j=1,\cdots,s \};
\end{eqnarray*}

\begin{eqnarray*}
R^{W}_l & = & \{2\alpha+2\sum\limits_{j=1}^b
a_{j}\delta_{1j}+4\sum\limits_{j=1}^s b_{j}\delta_{lj}
+\sum\limits_{j=1}^t c_{j}\delta'_{lj}\;|\;\alpha\in \dot{R}_s,
a_{j},b_j, c_j\in {\mathbf Z},
 \\
 & & \mbox{ and there is exactly one odd }c_j,j=1,\cdots,t \}.
\end{eqnarray*}
and $R^W_m$ is the following set:

$\mathrm{(4.1)}$\ for $l=2$,
\begin{eqnarray*}
R^{W}_m & = & \{\alpha+\sum\limits_{j=1}^b
a_{j}\delta_{1j}+2\sum\limits_{j=1}^s b_{j}\delta_{lj}
+\sum\limits_{j=1}^t c_{j}\delta'_{lj}\;|\;\alpha\in \dot{R}_l,
a_{j},b_j, c_j\in {\mathbf Z},
 \\
 & & \mbox{ and there is at most one odd }a_j,j=1\cdots,b \};
\end{eqnarray*}

$\mathrm{(4.2)}$\ for $l>2$,
$$R^{W}_m=\{\alpha+\sum\limits_{j=1}^b a_{j}\delta_{1j}+2\sum\limits_{j=1}^s b_{j}\delta_{lj}
+\sum\limits_{j=1}^t c_{j}\delta'_{lj}\;|\;\alpha\in \dot{R}_l,
a_{j},b_j, c_j\in {\mathbf Z}\}.$$

\end{theorem}

Before proving this theorem, we first show the following lemmas.

\begin{lemma}\label{form-of-Weyl-root}
If $\Delta\subseteq \dot{R}+rad(\Gamma)$, then
$R^{W}_s\subseteq\dot{R}_s+rad(\Gamma)$ and
$R^{W}_l\subseteq\dot{R}_l+rad(\Gamma)$.
\end{lemma}
\begin{proof}
For $\alpha, \beta\in\dot{R}+rad(\Gamma)$, let
$\alpha=\dot{\alpha}+\delta,\,\beta=\dot{\beta}+\delta'$, where
$\dot{\alpha}, \dot{\beta}\in\dot{R}$ and $\delta, \delta'\in
rad(\Gamma)$. Then
$s_\alpha(\beta)=\beta-(\alpha^{\vee},\,\beta)\alpha=\dot{\beta}+\delta'-
(\alpha^{\vee},\,\beta)(\dot{\alpha}+\delta)=\dot{\beta}+\delta'-
(\dot{\alpha}^{\vee},\,\dot{\beta})(\dot{\alpha}+\delta)=
s_{\dot{\alpha}}(\dot{\beta})+(\delta'-(\dot{\alpha}^{\vee},\,\dot{\beta})\delta)\in
\dot{R}+rad(\Gamma)$. Therefore if $\beta\in\dot{R}_s+rad(\Gamma)$,
then $s_\alpha(\beta)\in\dot{R}_s+rad(\Gamma)$; if
$\beta\in\dot{R}_l+rad(\Gamma)$, then
$s_\alpha(\beta)\in\dot{R}_l+rad(\Gamma)$. So by induction on the
number of reflections, we have
$R^{W}_s\subseteq\dot{R}_s+rad(\Gamma)$ and
$R^{W}_l\subseteq\dot{R}_l+rad(\Gamma)$.
\end{proof}

\begin{lemma}\label{plus-minus}
If $\alpha\in R^{W}$ and $\delta\in rad(\Gamma)$, then
$\alpha+\delta\in R^{W}$ if and only if $\alpha-\delta\in R^{W}$.
\end{lemma}
\begin{proof}
The result follows since $\alpha-\delta=-s_\alpha(\alpha+\delta)$.
\end{proof}

\begin{lemma}\label{mongline}
For $\alpha, \beta\in R^{W}$ with $(\alpha^\vee, \beta)=-1$ and
$\delta\in rad(\Gamma)$, if $\alpha+\delta\in R^{W}$,  then
$\beta+\delta\in R^{W}$. As a consequence, for any $\alpha\in
\dot{R}$, we have $\alpha+\delta_{1j}\in R^{W}$ for $j=1, \cdots,
b$.
\end{lemma}
\begin{proof}
Clearly $\beta+\delta=s_\alpha s_{\alpha+\delta}(\beta)\in R^{W}$.
Furthermore, note that $\alpha_1+\delta_{1j}=\alpha_{1j}\in R^{W}$
for $j=1, \cdots, b$. For any $\alpha_i, i=2, \cdots, l$, there is
an $\alpha_k, k<i$, such that $(\alpha_k^{\vee},\, \alpha_i)=-1$. By
induction, we have $\alpha_i+\delta_{1j}\in R^{W}$ for $i=1, \cdots,
l$ and $j=1, \cdots, b$. Therefore, for any $\alpha\in \dot{R}$, we
have $\alpha+\delta_{1j}\in R^{W}$ for $j=1, \cdots, b$.
\end{proof}

\begin{lemma}\label{even}
{\rm (1)}\ For $C_l(l\geq 2)$,  if $\alpha\in \dot{R}$ and
$\beta\in \dot{R_l}$, then $(\alpha^\vee,\,\beta)$ is even;\\
{\rm (2)}\ for $B_l(l\geq 3)$, if $\alpha\in \dot{R_s}$ and
$\beta\in \dot{R}$, then $(\alpha^\vee,\,\beta)$ is even;\\
{\rm (3)}\ for $F_4$, if $\alpha\in \dot{R_s}$ and $\beta\in
\dot{R_l}$, then $(\alpha^\vee,\,\beta)$ is even;\\
{\rm (4)}\ for $G_2$, if $\alpha\in \dot{R_s}$ and $\beta\in
\dot{R_l}$, then $(\alpha^\vee,\,\beta)$ is divided by 3.
\end{lemma}
\begin{proof}
In all cases, when $\alpha, \beta$  are simple roots, the result is
true. Note that
\begin{eqnarray*}
((s_{\gamma_1}(\alpha))^{\vee}, s_{\gamma_2}(\beta))&=&(\alpha^\vee,
\beta)-(\alpha^\vee, \gamma_1)(\gamma_1^\vee, \beta)- (\alpha^\vee,
\gamma_2)(\gamma_2^\vee, \beta)\\
& &+(\alpha^\vee, \gamma_1)(\gamma_1^\vee, \gamma_2)(\gamma_2^\vee,
\beta).
\end{eqnarray*}
 So by
induction on the number of reflections, we obtain the result.
\end{proof}

\vspace{0.4cm}
\begin{proof}[Proof of theorem~\ref{weyl-root}]
(1) Consider the type $X_l^{(d, 0, 0)}$, where $X_l=A_l(l\geq 2)$,
 $D_l$, $E_6$, $E_7$, $E_8$. Since
$\alpha_{1j}=\alpha_1+\delta_{1j}$ for $j=1, \cdots, s$, we have
$\alpha\in\dot{R}+rad(\Gamma)$ for any $\alpha\in\Delta$. By
Lemma~\ref{form-of-Weyl-root},
we obtain $R^{W}\subseteq\dot{R}+rad(\Gamma)$.\\
    \indent On the other hand, it suffices to show that
$\{\alpha_i+\sum\limits_{j=1}^d a_{j}\delta_{1j}\;|\;i=1, 2, \cdots,
l,\ a_{j}\in {\mathbf Z}\}\subseteq R^{W}$. For any $\alpha_i, i=2,
\cdots, l$, there is an $\alpha_k, k<i$, such that
$(\alpha_k^{\vee},\, \alpha_i)=-1$. So by Lemma~\ref{mongline}, we
only need to show that $S_1:=\{\alpha_1+\sum\limits_{j=1}^d
a_{j}\delta_{1j}\;|\;a_{j}\in {\mathbf Z}\}\subseteq R^{W}$. Given
any $\gamma\in S_1$, we have
$s_{\alpha_2}s_{\alpha_2+\delta_{1j}}(\gamma)=\gamma+\delta_{1j}$.
Therefore if $\gamma\in R^{W}$, we have $\gamma\pm\delta_{1j}\in
R^{W}$ for $j=1, \cdots, d$. By induction on the number
$n=\sum\limits_{j=1}^d |a_{j}|$, we obtain that $S_1\subseteq
R^{W}$.

(2) Consider the type $X_l^{(b, s, 0)}$, where $X_l=C_l(l\geq 3),
F_4, G_2$. Since $\alpha_{lj}=\alpha_l+\delta_{lj}$ for $j=1,
\cdots, s$, we have $\Delta_s\subseteq\dot{R}_s+rad(\Gamma)$. By
Lemma~\ref{form-of-Weyl-root},
we obtain $R^{W}_s\subseteq\dot{R}_s+rad(\Gamma)$.\\
    \indent On the other hand, by Lemma~\ref{mongline}, we only need to prove that
$S_2:=\{\alpha_l+\sum\limits_{j=1}^b
a_{j}\delta_{1j}+\sum\limits_{j=1}^s b_{j} \delta_{lj}\;|\;a_{j},
b_j\in {\mathbf Z}\}\subseteq R^{W}_s$. By Lemma~\ref{mongline}, we
have $\alpha_i+\delta_{1j}\in R^{W}$ for $i=1, \cdots,l,$ and $ j=1,
\cdots, b$. For $C_l(l\geq 3)$ and $F_4$, as
$\alpha_l+\delta_{lj}=\alpha_{lj}\in R^{W}$ for $j=1, \cdots, s$, by
Lemma~\ref{mongline}, we have $\alpha_{l-1}+\delta_{lj}\in R^{W}$
for $j=1, \cdots, s$; for $G_2$, as
$\alpha_2+\delta_{2j}=\alpha_{2j}\in R^{W}$ for $j=1, \cdots, s$, we
have
$\alpha_{1}+3\delta_{2j}=s_{\alpha_2}s_{\alpha_2+\delta_{2j}}(\alpha_1)\in
R^{W}$ for $j=1, \cdots, s$. Given any $\gamma\in S_2$, we have
$s_{\alpha_{l-1}}s_{\alpha_{l-1}+\delta_{1j}}(\gamma)=\gamma+\delta_{1j}$.
Furthermore, for $C_l(l\geq 3)$ and $F_4$, we have
$s_{\alpha_{l-1}}s_{\alpha_{l-1}+\delta_{lj}}(\gamma)=\gamma+\delta_{lj}$;
for $G_2$ we have
$s_{\alpha_{1}}s_{\alpha_{1}+3\delta_{2j}}s_{\alpha_{2}}s_{\alpha_{2}+\delta_{2j}}(\gamma)=\gamma+\delta_{2j}$.
So, if $\gamma\in R^{W}$, we have $\gamma\pm\delta_{1j},
\gamma\pm\delta_{2k}\in R^{W}$ for all $j=1, \cdots, b$, and $k=1,
\cdots, s$. By induction on the number $n=\sum\limits_{j=1}^b
|a_{j}|+\sum\limits_{j=1}^s |b_{j}|$, we obtain $S_2\subseteq
R^{W}_s$.

(2.1) $X_l=C_l(l\geq 3)$. Let $S_3$ be the set on the right-hand
side. Since $\alpha_{1j}=\alpha_1+\delta_{1j}$ for $j=1, \cdots, b$,
we have $\Delta_l\subseteq S_3$. Moreover, for any
$\beta=\alpha+\sum\limits_{j=1}^b
a_{j}\delta_{1j}+2\sum\limits_{j=1}^s b_{j}\delta_{lj}\in S_3$ with
$\alpha\in\dot{R}_l$, we have
$$s_{\alpha_{1k}}(\beta)=s_{\alpha_{1k}}(\alpha)+\sum\limits_{j=1}^b a_{j}\delta_{1j}+
2\sum\limits_{j=1}^s b_{j}\delta_{lj}=
 s_{\alpha_{1}}(\alpha)-(\alpha_{1}^\vee,
\alpha)\delta_{1k}+\sum\limits_{j=1}^b
a_{j}\delta_{1j}+2\sum\limits_{j=1}^s b_{j}\delta_{lj};$$
$$s_{\alpha_{lk}}(\beta)=s_{\alpha_{lk}}(\alpha)+\sum\limits_{j=1}^b a_{j}\delta_{1j}+
2\sum\limits_{j=1}^s b_{j}\delta_{lj}=
 s_{\alpha_{l}}(\alpha)-(\alpha_{l}^\vee,
\alpha)\delta_{lk}+\sum\limits_{j=1}^b
a_{j}\delta_{1j}+2\sum\limits_{j=1}^s b_{j}\delta_{lj}.$$ By
Lemma~\ref{even}, both $(\alpha_{1}^\vee, \alpha)$ and
$(\alpha_{l}^\vee, \alpha)$ are even. Then we have
$$s_{\alpha_{11}}(\beta), \cdots, s_{\alpha_{1b}}(\beta), s_{\alpha_{l1}}(\beta),
\cdots, s_{\alpha_{ls}}(\beta)\in S_3,$$ and
$s_{\alpha_{i}}(\beta)=s_{\alpha_{i}}(\alpha)+\sum\limits_{j=1}^b
a_{j}\delta_{1j}+ 2\sum\limits_{j=1}^s b_{j}\delta_{lj}\in S_3$ for
$i=1, 2, \cdots, l$.
So by induction on the number of reflections, we obtain $R^{W}_l\subseteq S_3$.\\
    \indent Next we prove $S_3\subseteq R^{W}_l$. It suffices to show that
$S_3'\subseteq R^{W}_l$, where $S_3':=\{\alpha_1+\sum\limits_{j=1}^b
a_{j}\delta_{1j}+2\sum\limits_{j=1}^s b_{j}\delta_{lj}\;|\; a_{j},
b_j\in {\mathbf Z},\mbox{ and there is at most one odd
}a_j,j=1,\cdots,b \}$. By Lemma~\ref{mongline}, we have
$\alpha_i+\delta_{1j}\in R^{W}$ for $i=1, \cdots, l$ and $j=1,
\cdots, b$. And since $\alpha_l+\delta_{lj}=\alpha_{lj}\in R^{W}$
for $j=1, \cdots, s$, by Lemma~\ref{mongline}, we have
$\alpha_i+\delta_{lj}\in R^{W}$ for $i=l-1, \cdots, 2 $ and $ j=1,
\cdots, s$. For any $\gamma\in S_3'$, we have
$s_{\alpha_{2}}s_{\alpha_{2}+\delta_{1j}}(\gamma)=\gamma+2\delta_{1j}$
and
$s_{\alpha_{2}}s_{\alpha_{2}+\delta_{lj}}(\gamma)=\gamma+2\delta_{lj}$.
Therefore, if $\gamma\in R^{W}$, we have $\gamma\pm2\delta_{1j},
\gamma\pm2\delta_{lk}\in R^{W}$ for all $j=1, \cdots, b $ and $ k=1,
\cdots, s$. By induction on the number $n=\sum\limits_{j=1}^b
|a_{j}|+2\sum\limits_{j=1}^s |b_{j}|$, we obtain $S_3'\subseteq
R^{W}_l$.

(2.2) $X_l=F_4$. Let $S_4$ be the set on the right-hand side. Since
$\alpha_{1j}=\alpha_1+\delta_{1j}$ for $j=1, \cdots, b$, we have
$\Delta_l\subseteq S_4$. Moreover, for any
$\beta=\alpha+\sum\limits_{j=1}^b
a_{j}\delta_{1j}+2\sum\limits_{j=1}^s b_{j}\delta_{4j}\in S_4$ with
$\alpha\in\dot{R}_l$, we have
$$s_{\alpha_{1k}}(\beta)=s_{\alpha_{1k}}(\alpha)+\sum\limits_{j=1}^b a_{j}\delta_{1j}+
2\sum\limits_{j=1}^s b_{j}\delta_{4j}=
 s_{\alpha_{1}}(\alpha)-(\alpha_{1}^\vee,
\alpha)\delta_{1k}+\sum\limits_{j=1}^b
a_{j}\delta_{1j}+2\sum\limits_{j=1}^s b_{j}\delta_{4j};$$
$$s_{\alpha_{4k}}(\beta)=s_{\alpha_{4k}}(\alpha)+\sum\limits_{j=1}^b a_{j}\delta_{1j}+
2\sum\limits_{j=1}^s b_{j}\delta_{4j}=
 s_{\alpha_{4}}(\alpha)-(\alpha_{4}^\vee,
\alpha)\delta_{4k}+\sum\limits_{j=1}^b
a_{j}\delta_{1j}+2\sum\limits_{j=1}^s b_{j}\delta_{4j}.$$ By
Lemma~\ref{even}, $(\alpha_{4}^\vee, \alpha)$ are even. Then we have
$$s_{\alpha_{11}}(\beta), \cdots, s_{\alpha_{1b}}(\beta), s_{\alpha_{41}}(\beta), \cdots,
s_{\alpha_{4s}}(\beta)\in S_4,$$ and
$s_{\alpha_{i}}(\beta)=s_{\alpha_{i}}(\alpha)+\sum\limits_{j=1}^b
a_{j}\delta_{1j}+ 2\sum\limits_{j=1}^s b_{j}\delta_{4j}\in S_4$ for
$i=1, 2, 3, 4$.
So by induction on the number of reflections, we obtain $R^{W}_l\subseteq S_4$.\\
    \indent It remains to prove $S_4\subseteq R^{W}_l$. By Lemma~\ref{mongline}, it suffices to show
$S_4':=\{\alpha_2+\sum\limits_{j=1}^b
a_{j}\delta_{1j}+2\sum\limits_{j=1}^s b_{j}\delta_{4j}\;|\; a_{j},
b_j\in {\mathbf Z}\}\subseteq R^{W}_l$. Since
$\alpha_4+\delta_{4j}=\alpha_{4j}\in R^{W}$ for $j=1, \cdots, s$, we
have $\alpha_3+\delta_{4j}\in R^{W}$ by Lemma~\ref{mongline}. For
any $\gamma\in S_4'$, we have
$s_{\alpha_{1}}s_{\alpha_{1j}}(\gamma)=\gamma+\delta_{1j}$ and
$s_{\alpha_{3}}s_{\alpha_{3}+\delta_{4j}}(\gamma)=\gamma+2\delta_{4j}$.
Thus, if $\gamma\in R^{W}$, we have $\gamma\pm\delta_{1j},
\gamma\pm2\delta_{4k}\in R^{W}$ for all $j=1, \cdots, b$ and $ k=1,
\cdots, s$. By induction on the number $n=\sum\limits_{j=1}^b
|a_{j}|+2\sum\limits_{j=1}^s |b_{j}|$, we have $S_4'\subseteq
R^{W}_l$.

(2.3) $X_l=G_2$. Let $S_5$ be the set on the right-hand side. Since
$\alpha_{1j}=\alpha_1+\delta_{1j}$ for $j=1, \cdots, b$, we obtain
$\Delta_l\subseteq S_5$. And for any
$\beta=\alpha+\sum\limits_{j=1}^b
a_{j}\delta_{1j}+3\sum\limits_{j=1}^s b_{j}\delta_{2j}\in S_5$,
where $\alpha\in\dot{R}_l$, we have
$$s_{\alpha_{1k}}(\beta)=s_{\alpha_{1k}}(\alpha)+\sum\limits_{j=1}^b a_{j}\delta_{1j}+
3\sum\limits_{j=1}^s b_{j}\delta_{2j}=
 s_{\alpha_{1}}(\alpha)-(\alpha_{1}^\vee,
\alpha)\delta_{1k}+\sum\limits_{j=1}^b
a_{j}\delta_{1j}+3\sum\limits_{j=1}^s b_{j}\delta_{2j};$$
$$s_{\alpha_{2k}}(\beta)=s_{\alpha_{2k}}(\alpha)+\sum\limits_{j=1}^b a_{j}\delta_{1j}+
3\sum\limits_{j=1}^s b_{j}\delta_{2j}=
 s_{\alpha_{2}}(\alpha)-(\alpha_{2}^\vee,
\alpha)\delta_{2k}+\sum\limits_{j=1}^b
a_{j}\delta_{1j}+3\sum\limits_{j=1}^s b_{j}\delta_{2j}.$$ By
Lemma~\ref{even}, $(\alpha_{2}^\vee, \alpha)$ is divided by $3$.
Then we obtain
$$s_{\alpha_{11}}(\beta), \cdots, s_{\alpha_{1b}}(\beta), s_{\alpha_{21}}(\beta), \cdots, s_{\alpha_{2s}}(\beta)\in S_5,$$
and
$s_{\alpha_{i}}(\beta)=s_{\alpha_{i}}(\alpha)+\sum\limits_{j=1}^b
a_{j}\delta_{1j}+ 3\sum\limits_{j=1}^s b_{j}\delta_{2j}\in S_5$ for
$i=1, 2$.
So by induction on the number of reflections, we get $R^{W}_l\subseteq S_5$.\\
    \indent Next we prove $S_5\subseteq R^{W}_l$. We only need to show that
$S_5'\subseteq R^{W}_l$, where $S_5':=\{\alpha_1+\sum\limits_{j=1}^b
a_{j}\delta_{1j}+3\sum\limits_{j=1}^s b_{j}\delta_{2j}\;|\; a_{j},
b_j\in {\mathbf Z}\}$. By Lemma~\ref{mongline}, we have
$\alpha_2+\delta_{1j}\in R^{W}$ for $j=1, \cdots, b$. For arbitrary
$\gamma\in S_5'$, we get
$s_{\alpha_{1}}s_{\alpha_{1j}}s_{\alpha_{2}}s_{\alpha_{2}+
\delta_{1j}}(\gamma)=\gamma+\delta_{1j}$ and
$s_{\alpha_{2}}s_{\alpha_{2j}}(\gamma)=\gamma+3\delta_{2j}$.
Therefore, if $\gamma\in R^{W}$, we obtain $\gamma\pm\delta_{1j},
\gamma\pm3\delta_{2k}\in R^{W}$ for all $j=1, \cdots, b $ and $ k=1,
\cdots, s$. By induction on the number $n=\sum\limits_{j=1}^b
|a_{j}|+3\sum\limits_{j=1}^s |b_{j}|$, we have $S_5'\subseteq
R^{W}_l$.

(3) Consider the type $A_1^{(b, 0, t)}$. Let
\begin{eqnarray*} S_6:& = &
\{\pm\alpha_1+\sum\limits_{j=1}^b
a_{j}\delta_{1j}+\sum\limits_{j=1}^t b_{j}
\delta'_{1j}\;|\;a_{j},b_j\in {\mathbf Z}, \\
 & & \mbox{ and there is at
most one odd }a_j,j=1\cdots,b \}.
\end{eqnarray*} Since
$\alpha_{1j}=\alpha_1+\delta_{1j}, j=1, \cdots, b$, we have
$\Delta_s\subseteq S_6$. And for any
$\beta=\alpha+\sum\limits_{j=1}^b
a_{j}\delta_{1j}+\sum\limits_{j=1}^t b_{j}\delta_{1j}'\in S_6$ with
$\alpha=\pm\alpha_1$, we obtain
$$s_{\alpha_{1k}}(\beta)=s_{\alpha_{1k}}(\alpha)+\sum\limits_{j=1}^b a_{j}\delta_{1j}+
\sum\limits_{j=1}^t b_{j}\delta_{1j}'=
s_{\alpha_{1}}(\alpha)-(\alpha_{1}^\vee,
\alpha)\delta_{1k}+\sum\limits_{j=1}^b
a_{j}\delta_{1j}+\sum\limits_{j=1}^t b_{j}\delta_{1j}';$$
$$s_{\beta_{1k}}(\beta)=s_{\beta_{1k}}(\alpha)+\sum\limits_{j=1}^b a_{j}\delta_{1j}+
\sum\limits_{j=1}^t b_{j}\delta_{1j}'=
 s_{\alpha_{1}}(\alpha)-\frac{1}{2}(\alpha_{1}^\vee,
\alpha)\delta_{1k}'+\sum\limits_{j=1}^b
a_{j}\delta_{1j}+\sum\limits_{j=1}^t b_{j}\delta_{1j}'.$$ Note that
$(\alpha_{1}^\vee, \alpha)=\pm2$. Then we get
$$s_{\alpha_{11}}(\beta), \cdots, s_{\alpha_{1b}}(\beta), s_{\beta'_{11}}(\beta), \cdots,
s_{\beta'_{1t}}(\beta)\in S_6,$$ and
$s_{\alpha_{1}}(\beta)=s_{\alpha_{1}}(\alpha)+\sum\limits_{j=1}^b
a_{j}\delta_{1j}+ \sum\limits_{j=1}^t b_{j}\delta_{1j}'\in S_6$.
So by induction on the number of reflections, we obtain $R^{W}_s\subseteq S_6$.\\
    \indent It remains to show $S_6\subseteq R^{W}_s$. We only need to show
$S_6'\subseteq R^{W}_s$, where $S_6':=\{\alpha_1+\sum\limits_{j=1}^b
a_{j}\delta_{1j}+\sum\limits_{j=1}^t
b_{j}\delta_{1j}'\;|\;a_{j},b_j\in {\mathbf Z}, \mbox{ and there is
at most one odd }a_j,j=1\cdots,b\}$. For arbitrary $\gamma\in S_6'$,
we obtain
$s_{\alpha_{1j}}s_{\alpha_{1}}(\gamma)=\gamma+2\delta_{1j}$ and
$s_{\beta_{1j}}s_{\alpha_{1}}(\gamma)=\gamma+\delta_{1j}'$.
Therefore, if $\gamma\in R^{W}$, we have $\gamma\pm2\delta_{1j},
\gamma\pm\delta_{1k}'\in R^{W}$ for all $j=1, \cdots, b $ and $ k=1,
\cdots, s$. Moreover, $\alpha_1+\delta_{1j}=\alpha_{1j}\in R^{W}_s$
for $j=1, \cdots, b$. By induction on the number
$n=\sum\limits_{j=1}^b |a_{j}|+\sum\limits_{j=1}^t |b_{j}|$, we get
$S_6'\subseteq R^{W}_s$.

So $R^{W}_s=S_6=\{\pm\alpha_1+\sum\limits_{j=1}^b
a_{j}\delta_{1j}+\sum\limits_{j=1}^t b_{j}
\delta'_{1j}\;|\;a_{j},b_j\in {\mathbf Z},  \mbox{ and there is at
most one}\\
\mbox{odd }a_j,j=1\cdots,b \}$.

Let $S_7:=\{\pm2\alpha_1+4\sum\limits_{j=1}^b
a_{j}\delta_{1j}+\sum\limits_{j=1}^t b_{j}\delta'_{1j}\;|\;
a_{j},b_j\in {\mathbf Z},  \mbox{ and there is exactly one odd}\\
b_j,j=1,\cdots,t \}$. Since $\beta_{1j}=2\alpha_1+\delta'_{1j}, j=1,
\cdots, t$, we have $\Delta_l\subseteq S_7$. And for any
$\beta=\alpha+4\sum\limits_{j=1}^b
a_{j}\delta_{1j}+\sum\limits_{j=1}^t b_{j}\delta'_{1j}\in S_7$ with
$\alpha=\pm2\alpha_1$, we have
$$s_{\alpha_{1k}}(\beta)=s_{\alpha_{1k}}(\alpha)+4\sum\limits_{j=1}^b a_{j}\delta_{1j}+
\sum\limits_{j=1}^t b_{j}\delta'_{1j}=
 s_{\alpha_{1}}(\alpha)-(\alpha_{1}^\vee,
\alpha)\delta_{1k}+4\sum\limits_{j=1}^b
a_{j}\delta_{1j}+\sum\limits_{j=1}^t b_{j}\delta'_{1j};$$
$$s_{\beta_{1k}}(\beta)=s_{\beta_{1k}}(\alpha)+4\sum\limits_{j=1}^b a_{j}\delta_{1j}+
\sum\limits_{j=1}^t b_{j}\delta'_{1j}=
s_{\alpha_{1}}(\alpha)-\frac{1}{2}(\alpha_{1}^\vee,
\alpha)\delta'_{1k}+4\sum\limits_{j=1}^b
a_{j}\delta_{1j}+\sum\limits_{j=1}^t b_{j}\delta'_{1j}.$$ Moreover,
$s_{\alpha_{1}}(\alpha)=-\alpha, \ (\alpha_{1}^\vee, \alpha)=\pm4$.
Then we obtain
$$s_{\alpha_{11}}(\beta), \cdots, s_{\alpha_{1b}}(\beta), s_{\beta_{11}}(\beta), \cdots,
s_{\beta_{1t}}(\beta)\in S_7,$$ and
$s_{\alpha_{1}}(\beta)=s_{\alpha_{1}}(\alpha)+4\sum\limits_{j=1}^b
a_{j}\delta_{1j}+ \sum\limits_{j=1}^t b_{j}\delta'_{1j}\in S_7$.
So by induction on the number of reflections, we get $R^{W}_l\subseteq S_7$.\\
    \indent Next we prove $S_7\subseteq R^{W}_l$. It suffices to show
$S_7'\subseteq R^{W}_l$, where
$S_7':=\{2\alpha_1+4\sum\limits_{j=1}^b
a_{j}\delta_{1j}+\sum\limits_{j=1}^t b_{j}\delta'_{1j}\;|\;a_{j},b_j
\in {\mathbf Z},\mbox{ and there is exactly one odd }
b_j,j=1,\cdots,t\}$. For any $\gamma\in S_7'$, we get
$s_{\alpha_{1j}}s_{\alpha_{1}}(\gamma)=\gamma+4\delta_{1j}$ and
$s_{\beta_{1j}}s_{\alpha_{1}}(\gamma)=\gamma+2\delta'_{1j}$. So, if
$\gamma\in R^{W}$, then $\gamma\pm4\delta_{1j},
\gamma\pm2\delta'_{1k}\in R^{W}$ for all $j=1, \cdots, b $ and $
k=1, \cdots, s$. Notice that $2\alpha_1+\delta'_{1j}=\beta_{1j}\in
R^{W}_l$ for $j=1, \cdots, t$. So by induction on the number
$n=4\sum\limits_{j=1}^b |a_{j}|+\sum\limits_{j=1}^t |b_{j}|$, we
obtain that $S_7'\subseteq R^{W}_l$.

Hence $R^{W}_l=S_7=\{\pm2\alpha_1+4\sum\limits_{j=1}^b
a_{j}\delta_{1j}+\sum\limits_{j=1}^t b_{j}\delta'_{1j}\;|\;
a_{j},b_j\in {\mathbf Z},  \mbox{ and there is exactly }\\
\mbox{one odd } b_j,j=1,\cdots,t \}$.

(4) Consider the type $B_l^{(b, s, t)}$. Let
$S_8:=\{\alpha+\sum\limits_{j=1}^b
a_{j}\delta_{1j}+\sum\limits_{j=1}^s b_{j}\delta_{lj}
+\sum\limits_{j=1}^t c_{j}\delta'_{lj}\;|\;\alpha\in \dot{R}_s,
a_{j},b_j, c_j\in {\mathbf Z}, \mbox{and there is at most one odd} \
b_j,j=1,\cdots,s \}$. Since $\alpha_{lj}=\alpha_l+\delta_{lj}, j=1,
\cdots, s$, we have $\Delta_s\subseteq S_8$. And for any
$\beta=\alpha+\sum\limits_{j=1}^b
a_{j}\delta_{1j}+\sum\limits_{j=1}^s b_{j}\delta_{lj}
+\sum\limits_{j=1}^t c_{j}\delta'_{lj}\in S_8$ with $\alpha\in
\dot{R}_s$, we obtain
\begin{eqnarray*}
s_{\alpha_{1k}}(\beta) & = &
s_{\alpha_{1k}}(\alpha)+\sum\limits_{j=1}^b a_{j}\delta_{1j}+
\sum\limits_{j=1}^s b_{j}\delta_{lj}+\sum\limits_{j=1}^t
c_{j}\delta'_{lj}\\
 & = &
 s_{\alpha_{1}}(\alpha)-(\alpha_{1}^\vee,
\alpha)\delta_{1k}+\sum\limits_{j=1}^b
a_{j}\delta_{1j}+\sum\limits_{j=1}^s b_{j}\delta_{lj}
+\sum\limits_{j=1}^t c_{j}\delta'_{lj};
\end{eqnarray*}
\begin{eqnarray*}
 s_{\alpha_{lk}}(\beta) & = &
s_{\alpha_{lk}}(\alpha)+\sum\limits_{j=1}^b a_{j}\delta_{1j}+
\sum\limits_{j=1}^s b_{j}\delta_{lj}+\sum\limits_{j=1}^t
c_{j}\delta'_{lj}\\
 & = &
 s_{\alpha_{l}}(\alpha)-(\alpha_{l}^\vee,
\alpha)\delta_{lk}+\sum\limits_{j=1}^b
a_{j}\delta_{1j}+\sum\limits_{j=1}^s b_{j}\delta_{lj}
+\sum\limits_{j=1}^t c_{j}\delta'_{lj};
\end{eqnarray*}
\begin{eqnarray*}
 s_{\beta_{lk}}(\beta) & = &
s_{\beta_{lk}}(\alpha)+\sum\limits_{j=1}^b a_{j}\delta_{1j}+
\sum\limits_{j=1}^s b_{j}\delta_{lj}+\sum\limits_{j=1}^t
c_{j}\delta'_{lj}\\
 & = &
 s_{\alpha_{l}}(\alpha)-\frac{1}{2}(\alpha_{l}^\vee,
\alpha)\delta'_{lk}+\sum\limits_{j=1}^b
a_{j}\delta_{1j}+\sum\limits_{j=1}^s b_{j}\delta_{lj}
+\sum\limits_{j=1}^t c_{j}\delta'_{lj}.
\end{eqnarray*}
By Lemma~\ref{even}, $(\alpha_{l}^\vee, \alpha)$ is even. Thus we
have
$$s_{\alpha_{11}}(\beta), \cdots, s_{\alpha_{1b}}(\beta), s_{\alpha_{l1}}(\beta), \cdots, s_{\alpha_{ls}}(\beta),
s_{\beta_{l1}}(\beta), \cdots, s_{\beta_{lt}}(\beta)\in S_8,$$ and
$s_{\alpha_{i}}(\beta)=s_{\alpha_{i}}(\alpha)+\sum\limits_{j=1}^b
a_{j}\delta_{1j}+\sum\limits_{j=1}^s b_{j}\delta_{lj}
+\sum\limits_{j=1}^t c_{j}\delta'_{lj}\in S_8$ for $i=1, \cdots, l$.
So by induction on the number of reflections, we get $R^{W}_s\subseteq S_8$.\\
    \indent It remains to prove $S_8\subseteq R^{W}_s$. We only need to show that
$S_8'\subseteq R^{W}_s$, where $S_8':=\{\alpha_l+\sum\limits_{j=1}^b
a_{j}\delta_{1j}+\sum\limits_{j=1}^s b_{j}\delta_{lj}
+\sum\limits_{j=1}^t c_{j}\delta'_{lj}\;|\;a_{j},b_j, c_j\in
{\mathbf Z}, \mbox{ and there is at most one odd
}b_j,\\j=1,\cdots,s\}$. By Lemma~\ref{mongline}, we have
$\alpha_i+\delta_{1j}\in R^{W}$ for $i=1, \cdots, l, j=1, \cdots,
b$. For any $\gamma\in S_8'$, we get
$s_{\alpha_{l-1}}s_{\alpha_{l-1}+\delta_{1j}}(\gamma)=\gamma+\delta_{1j}$,
$s_{\alpha_{lj}}s_{\alpha_{l}}(\gamma)=\gamma+2\delta_{lj}$ and
$s_{\beta_{lj}}s_{\alpha_{l}}(\gamma)=\gamma+\delta'_{lj}$.
Therefore, if $\gamma\in R^{W}$, we have $\gamma\pm\delta_{1i},
\gamma\pm2\delta_{lj}, \gamma\pm\delta'_{lk}\in R^{W}$ for all $i=1,
\cdots, b, j=1, \cdots, s $ and $ k=1, \cdots, t$. Note that
$\alpha_l+\delta_{lj}=\alpha_{lj}\in R^{W}_s$ for $j=1, \cdots, s$.
By induction on the number $n=\sum\limits_{j=1}^b
|a_{j}|+\sum\limits_{j=1}^s |b_{j}|+\sum\limits_{j=1}^t |c_{j}|$, we
obtain $S_8'\subseteq R^{W}_s$.

So $R^{W}_s=S_8=\{\alpha+\sum\limits_{j=1}^b
a_{j}\delta_{1j}+\sum\limits_{j=1}^s b_{j}\delta_{lj}
+\sum\limits_{j=1}^t c_{j}\delta'_{lj}\;|\;\alpha\in \dot{R}_s,
a_{j},b_j, c_j\in {\mathbf Z}, \mbox{and there
 is at most one odd} \ b_j,j=1,\cdots,s \}$.

Let $S_9:=\{2\alpha+2\sum\limits_{j=1}^b
a_{j}\delta_{1j}+4\sum\limits_{j=1}^s b_{j}\delta_{lj}
+\sum\limits_{j=1}^t c_{j}\delta'_{lj}\;|\;\alpha\in \dot{R}_s,
a_{j},b_j, c_j\in {\mathbf Z}, \mbox{and there }
 \\
\mbox{is exactly one odd }c_j,j=1,\cdots,t \}$. Since
$\beta_{lj}=2\alpha_l+\delta'_{lj}, j=1, \cdots, t$, we get
$\Delta_l\subseteq S_9$. And for any
$\beta=2\alpha+2\sum\limits_{j=1}^b
a_{j}\delta_{1j}+4\sum\limits_{j=1}^s b_{j}\delta_{lj}
+\sum\limits_{j=1}^t c_{j}\delta'_{lj}\in S_9$ with $\alpha\in
\dot{R}_s$, we have
\begin{eqnarray*}
s_{\alpha_{1k}}(\beta) & = &
2s_{\alpha_{1k}}(\alpha)+2\sum\limits_{j=1}^b a_{j}\delta_{1j}+
4\sum\limits_{j=1}^s b_{j}\delta_{lj}+\sum\limits_{j=1}^t c_{j}\delta'_{lj}\\
 & = & 2s_{\alpha_{1}}(\alpha)-2(\alpha_{1}^\vee,
\alpha)\delta_{1k}+2\sum\limits_{j=1}^b
a_{j}\delta_{1j}+4\sum\limits_{j=1}^s b_{j}\delta_{lj}
+\sum\limits_{j=1}^t c_{j}\delta'_{lj};
\end{eqnarray*}
\begin{eqnarray*}
s_{\alpha_{lk}}(\beta) & = &
2s_{\alpha_{lk}}(\alpha)+2\sum\limits_{j=1}^b a_{j}\delta_{1j}+
4\sum\limits_{j=1}^s b_{j}\delta_{lj}+\sum\limits_{j=1}^t c_{j}\delta'_{lj}\\
 & = & 2s_{\alpha_{l}}(\alpha)-2(\alpha_{l}^\vee,
\alpha)\delta_{lk}+2\sum\limits_{j=1}^b
a_{j}\delta_{1j}+4\sum\limits_{j=1}^s b_{j}\delta_{lj}
+\sum\limits_{j=1}^t c_{j}\delta'_{lj};
\end{eqnarray*}
\begin{eqnarray*}
 s_{\beta_{lk}}(\beta) & = &
2s_{\beta_{lk}}(\alpha)+2\sum\limits_{j=1}^b a_{j}\delta_{1j}+
4\sum\limits_{j=1}^s b_{j}\delta_{lj}+\sum\limits_{j=1}^t c_{j}\delta'_{lj}\\
 & = & 2s_{\alpha_{l}}(\alpha)-(\alpha_{l}^\vee,
\alpha)\delta'_{lk}+2\sum\limits_{j=1}^b
a_{j}\delta_{1j}+4\sum\limits_{j=1}^s b_{j}\delta_{lj}
+\sum\limits_{j=1}^t c_{j}\delta'_{lj}.
\end{eqnarray*}
By Lemma~\ref{even}, $(\alpha_{l}^\vee, \alpha)$ is even. Then we
obtain
$$s_{\alpha_{11}}(\beta), \cdots, s_{\alpha_{1b}}(\beta), s_{\alpha_{l1}}(\beta), \cdots, s_{\alpha_{ls}}(\beta),
s_{\beta_{l1}}(\beta), \cdots, s_{\beta_{lt}}(\beta)\in S_9,$$ and
$s_{\alpha_{i}}(\beta)=2s_{\alpha_{i}}(\alpha)+2\sum\limits_{j=1}^b
a_{j}\delta_{1j}+4\sum\limits_{j=1}^s b_{j}\delta_{lj}
+\sum\limits_{j=1}^t c_{j}\delta'_{lj}\in S_9$ for $i=1, \cdots, l$.
So by induction on the number of reflections, we obtain $R^{W}_l\subseteq S_9$.\\
    \indent Next we prove $S_9\subseteq R^{W}_l$. It suffices to show that
$S_9'\subseteq R^{W}_l$, where
$S_9':=\{2\alpha_l+2\sum\limits_{j=1}^b
a_{j}\delta_{1j}+4\sum\limits_{j=1}^s b_{j}\delta_{lj}
+\sum\limits_{j=1}^t c_{j}\delta'_{lj}\;|\;a_{j},b_j, c_j\in
{\mathbf Z}, \mbox{ and there is exactly one odd
}c_j,j=1,\cdots,t\}$. By Lemma~\ref{mongline}, we have
$\alpha_i+\delta_{1j}\in R^{W}$ for $i=1, \cdots, l$ and $j=1,
\cdots, b$. For any $\gamma\in S_9'$, we get
$s_{\alpha_{l-1}}s_{\alpha_{l-1}+
\delta_{1j}}(\gamma)=\gamma+2\delta_{1j}$,
$s_{\alpha_{lj}}s_{\alpha_{l}}(\gamma)=\gamma+4\delta_{lj}$ and
$s_{\beta_{lj}}s_{\alpha_{l}}(\gamma)=\gamma+2\delta'_{lj}$.
Therefore, if $\gamma\in R^{W}$, we obtain $\gamma\pm2\delta_{1i},
\gamma\pm4\delta_{lj}, \gamma\pm2\delta'_{lk}\in R^{W}$ for all
$i=1, \cdots, b, j=1, \cdots, s $ and $ k=1, \cdots, t$. Note that
$2\alpha_l+\delta'_{lj}=\beta_{lj}\in R^{W}_l$ for $j=1, \cdots, s$.
By induction on the number $n=2\sum\limits_{j=1}^b
|a_{j}|+4\sum\limits_{j=1}^s |b_{j}|+\sum\limits_{j=1}^t |c_{j}|$,
we have $S_9'\subseteq R^{W}_l$.

So $R^{W}_l=S_9=\{2\alpha+2\sum\limits_{j=1}^b
a_{j}\delta_{1j}+4\sum\limits_{j=1}^s b_{j}\delta_{lj}
+\sum\limits_{j=1}^t c_{j}\delta'_{lj}\;|\;\alpha\in \dot{R}_s,
a_{j},b_j, c_j\in {\mathbf Z}, \mbox{and there is exactly one odd
}c_j,j=1,\cdots,t \}$.

(4.1) $X_l=B_l(l=2)$. Let $S_{10}$ be the set on the right-hand
side. Since $\alpha_{1j}=\alpha_1+\delta_{1j}, j=1, \cdots, b$, we
have $\Delta_m\subseteq S_{10}$. And for any
$\beta=\alpha+\sum\limits_{j=1}^b
a_{j}\delta_{1j}+2\sum\limits_{j=1}^s b_{j}\delta_{2j}
+\sum\limits_{j=1}^t c_{j}\delta'_{2j}\in S_{10}$ with $\alpha\in
\dot{R}_l$, we have
\begin{eqnarray*}
s_{\alpha_{1k}}(\beta) & = &
s_{\alpha_{1k}}(\alpha)+\sum\limits_{j=1}^b a_{j}\delta_{1j}+
2\sum\limits_{j=1}^s b_{j}\delta_{2j} +\sum\limits_{j=1}^t
c_{j}\delta'_{2j}\\
 & = &
 s_{\alpha_{1}}(\alpha)-(\alpha_{1}^\vee,
\alpha)\delta_{1k}+\sum\limits_{j=1}^b
a_{j}\delta_{1j}+2\sum\limits_{j=1}^s b_{j}\delta_{2j}
+\sum\limits_{j=1}^t c_{j}\delta'_{2j};
\end{eqnarray*}
\begin{eqnarray*}
s_{\alpha_{2k}}(\beta) & = &
s_{\alpha_{2k}}(\alpha)+\sum\limits_{j=1}^b a_{j}\delta_{1j}+
2\sum\limits_{j=1}^s b_{j}\delta_{2j} +\sum\limits_{j=1}^t
c_{j}\delta'_{2j}\\
 & = &
 s_{\alpha_{2}}(\alpha)-(\alpha_{2}^\vee,
\alpha)\delta_{2k}+\sum\limits_{j=1}^b
a_{j}\delta_{1j}+2\sum\limits_{j=1}^s b_{j}\delta_{2j}
+\sum\limits_{j=1}^t c_{j}\delta'_{2j};
\end{eqnarray*}
\begin{eqnarray*}
s_{\beta_{2k}}(\beta) & = &
s_{\beta_{2k}}(\alpha)+\sum\limits_{j=1}^b a_{j}\delta_{1j}+
2\sum\limits_{j=1}^s b_{j}\delta_{2j} +\sum\limits_{j=1}^t
c_{j}\delta'_{2j}\\
 & = &
 s_{\alpha_{2}}(\alpha)-\frac{1}{2}(\alpha_{2}^\vee,
\alpha)\delta'_{2k}+\sum\limits_{j=1}^b
a_{j}\delta_{1j}+2\sum\limits_{j=1}^s b_{j}\delta_{2j}
+\sum\limits_{j=1}^t c_{j}\delta'_{2j}.
\end{eqnarray*}
By Lemma~\ref{even}, $(\alpha_{i}^\vee, \alpha)$ is even for $ i=1,
2$. Then we have
$$s_{\alpha_{11}}(\beta), \cdots, s_{\alpha_{1b}}(\beta), s_{\alpha_{21}}(\beta), \cdots, s_{\alpha_{2s}}(\beta),
s_{\beta_{21}}(\beta), \cdots, s_{\beta_{2t}}(\beta)\in S_{10},$$
and
$s_{\alpha_{i}}(\beta)=s_{\alpha_{i}}(\alpha)+\sum\limits_{j=1}^b
a_{j}\delta_{1j}+2\sum\limits_{j=1}^s b_{j}\delta_{2j}
+\sum\limits_{j=1}^t c_{j}\delta'_{2j}\in S_{10}$ for $i=1, 2$.
So by induction on the number of reflections, we obtain $R^{W}_m\subseteq S_{10}$.\\
    \indent Next we prove $S_{10}\subseteq R^{W}_m$. It suffices to show that
$S_{10}'\subseteq R^{W}_m$, where
$S_{10}':=\{\alpha_1+\sum\limits_{j=1}^b
a_{j}\delta_{1j}+2\sum\limits_{j=1}^s b_{j}\delta_{2j}
+\sum\limits_{j=1}^t c_{j}\delta'_{2j}\;|\;a_{j},b_j, c_j\in
{\mathbf Z}, \mbox{ and there is at most one odd }a_j,j=1\cdots,b
\}$. For any $\gamma\in S_{10}'$, we get
$s_{\alpha_{1j}}s_{\alpha_{1}}(\gamma)=\gamma+2\delta_{1j}$,
$s_{\alpha_{2}}s_{\alpha_{2j}}(\gamma)=\gamma+2\delta_{2j}$ and
$s_{\alpha_{2}}s_{\beta_{2j}}(\gamma)=\gamma+\delta'_{2j}$.
Therefore, if $\gamma\in R^{W}$,
 we have $\gamma\pm2\delta_{1i}, \gamma\pm2\delta_{2j}, \gamma\pm\delta'_{2k}\in R^{W}$
 for all $i=1, \cdots, b, j=1, \cdots, s$ and $k=1, \cdots, t$.
Note that $\alpha_1+\delta_{1j}=\alpha_{1j}\in R^{W}_m$ for $j=1,
\cdots, b$. Hence, by induction on the number $n=\sum\limits_{j=1}^b
|a_{j}|+2\sum\limits_{j=1}^s |b_{j}|+ \sum\limits_{j=1}^t |c_{j}|$,
we obtain $S_{10}'\subseteq R^{W}_m$.

(4.2) $X_l=B_l (l>2)$. Let $S_{11}$ be the set on the right-hand
side. Since $\alpha_{1j}=\alpha_1+\delta_{1j}, j=1, \cdots, b$, we
have $\Delta_m\subseteq S_{11}$. And for any
$\beta=\alpha+\sum\limits_{j=1}^b
a_{j}\delta_{1j}+2\sum\limits_{j=1}^s b_{j}\delta_{lj}
+\sum\limits_{j=1}^t c_{j}\delta'_{lj}\in S_{11}$ with $\alpha\in
\dot{R}_l$, we have
\begin{eqnarray*}
s_{\alpha_{1k}}(\beta) & = &
s_{\alpha_{1k}}(\alpha)+\sum\limits_{j=1}^b a_{j}\delta_{1j}+
2\sum\limits_{j=1}^s b_{j}\delta_{lj}+\sum\limits_{j=1}^t
c_{j}\delta'_{lj}\\
 & = &
 s_{\alpha_{1}}(\alpha)-(\alpha_{1}^\vee,
\alpha)\delta_{1k}+\sum\limits_{j=1}^b
a_{j}\delta_{1j}+2\sum\limits_{j=1}^s b_{j}\delta_{lj}
+\sum\limits_{j=1}^t c_{j}\delta'_{lj};
\end{eqnarray*}
\begin{eqnarray*}
s_{\alpha_{lk}}(\beta) & = &
s_{\alpha_{lk}}(\alpha)+\sum\limits_{j=1}^b a_{j}\delta_{1j}+
2\sum\limits_{j=1}^s b_{j}\delta_{lj} +\sum\limits_{j=1}^t
c_{j}\delta'_{lj}\\
 & = & s_{\alpha_{l}}(\alpha)-(\alpha_{l}^\vee,
\alpha)\delta_{lk}+\sum\limits_{j=1}^b
a_{j}\delta_{1j}+2\sum\limits_{j=1}^s b_{j}\delta_{lj}
+\sum\limits_{j=1}^t c_{j}\delta'_{lj};
\end{eqnarray*}
\begin{eqnarray*}
s_{\beta_{lk}}(\beta) & = &
s_{\beta_{lk}}(\alpha)+\sum\limits_{j=1}^b a_{j}\delta_{1j}+
2\sum\limits_{j=1}^s b_{j}\delta_{lj} +\sum\limits_{j=1}^t
c_{j}\delta'_{lj}\\
 & = & s_{\alpha_{l}}(\alpha)-\frac{1}{2}(\alpha_{l}^\vee,
\alpha)\delta'_{lk}+\sum\limits_{j=1}^b
a_{j}\delta_{1j}+2\sum\limits_{j=1}^s b_{j}\delta_{lj}
+\sum\limits_{j=1}^t c_{j}\delta'_{lj}.
\end{eqnarray*}
By Lemma~\ref{even}, $(\alpha_{l}^\vee, \alpha)$ is even. Then we
obtain
$$s_{\alpha_{11}}(\beta), \cdots, s_{\alpha_{1b}}(\beta), s_{\alpha_{l1}}(\beta), \cdots, s_{\alpha_{ls}}(\beta),
s_{\beta_{l1}}(\beta), \cdots, s_{\beta_{lt}}(\beta)\in S_{11},$$
and
$s_{\alpha_{i}}(\beta)=s_{\alpha_{i}}(\alpha)+\sum\limits_{j=1}^b
a_{j}\delta_{1j}+2\sum\limits_{j=1}^s b_{j}\delta_{lj}
+\sum\limits_{j=1}^t c_{j}\delta'_{lj}\in S_{11}$ for $i=1, \cdots,
l$.
So by induction on the number of reflections, we obtain $R^{W}_m\subseteq S_{11}$.\\
    \indent It remains to prove $S_{11}\subseteq R^{W}_m$. By Lemma~\ref{mongline}, we only need to show
$S_{11}'\subseteq R^{W}_m$, where
$S_{11}':=\{\alpha_{l-1}+\sum\limits_{j=1}^b
a_{j}\delta_{1j}+2\sum\limits_{j=1}^s b_{j}\delta_{lj}
+\sum\limits_{j=1}^t c_{j}\delta'_{lj}\;|\;a_{j},b_j, c_j\in
{\mathbf Z}\}$. By Lemma~\ref{mongline}, we have
$\alpha_i+\delta_{1j}\in R^{W}$ for $i=1, \cdots, l, j=1, \cdots,
b$. For any $\gamma\in S_{11}'$, we obtain
$s_{\alpha_{l-2}}s_{\alpha_{l-2}+\delta_{1j}}(\gamma)=\gamma+\delta_{1j}$,
$s_{\alpha_{l}}s_{\alpha_{lj}}(\gamma)=\gamma+2\delta_{lj}$ and
$s_{\alpha_{l}}s_{\beta_{lj}}(\gamma)=\gamma+\delta'_{lj}$.
Therefore, if $\gamma\in R^{W}$,
 we have $\gamma\pm\delta_{1i}, \gamma\pm2\delta_{lj}, \gamma\pm\delta'_{lk}\in R^{W}$ for all
 $i=1, \cdots, b, j=1, \cdots, s$ and $k=1, \cdots, t$.
By induction on the number $n=\sum\limits_{j=1}^b
|a_{j}|+2\sum\limits_{j=1}^s |b_{j}|+\sum\limits_{j=1}^t |c_{j}|$,
we obtain $S_{11}'\subseteq R^{W}_m$.
\end{proof}

\end{document}